\documentclass{amsart}

\usepackage{a4wide}
\usepackage{amsmath}
\usepackage{amsfonts}
\usepackage{rotating}
\usepackage{bbm}
\usepackage[all,cmtip]{xy}

\usepackage{amssymb}
\usepackage{graphicx}
\usepackage{dsfont}
\usepackage{color}
\usepackage{hyperref}
\usepackage{verbatim}
\usepackage{mathtools}
\usepackage{enumitem}
\usepackage{color}
\usepackage{epstopdf}
\usepackage{booktabs}

\usepackage{graphicx}
\usepackage[utf8]{inputenc}
\usepackage{kotex}
\usepackage{amsfonts}
\usepackage{amsmath}
\usepackage{amsthm}
\usepackage{amssymb}
\usepackage{tikz}
\usepackage[toc,page]{appendix}
\usetikzlibrary{patterns}
\usetikzlibrary{decorations.pathreplacing}
\usetikzlibrary{cd}
\usepackage[margin=1.0in]{geometry}
\usepackage{caption}
\usepackage{float}
\usepackage{mathtools}
\usepackage{pgfplots}
\usepackage{tkz-euclide}
\newcommand{\Z}{\mathbb{Z}}
\newcommand{\R}{\mathbb{R}}

\DeclareMathOperator{\lk}{lk}

\newtheorem{thm}{Theorem}[section]
\newtheorem{protothm}[thm]{Proto-Theorem}

\newtheorem{observation}[thm]{Observation}
\newtheorem{lemma}[thm]{Lemma}
\newtheorem{proposition}[thm]{Proposition}
\newtheorem{cor}[thm]{Corollary}
\newtheorem{defn}[thm]{Definition}
\newtheorem{remark}[thm]{Remark}
\newtheorem{corollary}[thm]{Corollary}

\newcommand{\plota}[1][]{
\begin{tikzpicture}
\draw [fill=lightgray] (2.5,0) circle (2);
\draw [fill=white] (2.5,0) circle (1.5);
\draw [fill=lightgray] (2.5,0) circle (1);
\draw [fill=white] (2.5,0) circle (0.5);

\draw [draw=none, fill=lightgray] (3.2,-0.25) rectangle (4.2,0.25);
\draw (3.46,-0.25) -- (3.985,-0.25);
\draw (3.46,0.25) -- (3.985,0.25);

\draw (2.6,0.2) node {$U_2$};
\draw (4.5,2) node {$U_1\#L_2$};
\draw (4,1.8) -- (3.45,1.17);
\draw (4.7,1.8) -- (3.3,0.6);
\end{tikzpicture}
}

\newcommand*{\halfway}{0.5*\pgfdecoratedpathlength+.5*2mm}

\newcommand{\plotb}[1][]{
\begin{tikzpicture}
\draw [fill=lightgray] (2.5,-5) circle (1);
\draw [fill=white] (2.5,-5) circle (0.75);
\draw [fill=lightgray] (2.5,-5) circle (0.5);
\draw [fill=white] (2.5,-5) circle (0.25);
\draw [draw=none, fill=lightgray] (2.8,-4.85) rectangle (3.35,-5.15);
\draw (2.98,-4.85) -- (3.25,-4.85);
\draw (2.98,-5.15) -- (3.25,-5.15);

\draw [fill=lightgray] (0,-2.5) circle (1);
\draw [fill=white] (0,-2.5) circle (0.75);
\draw [fill=lightgray] (0,-2.5) circle (0.5);
\draw [fill=white] (0,-2.5) circle (0.25);
\draw [draw=none, fill=lightgray] (-0.15,-2.2) rectangle (0.15,-1.65);
\draw (-0.15,-1.76) -- (-0.15,-2.03);
\draw (0.15,-1.76) -- (0.15,-2.03);

\draw [fill=lightgray] (-2.5,-5) circle (1);
\draw [fill=white] (-2.5,-5) circle (0.75);
\draw [fill=lightgray] (-2.5,-5) circle (0.5);
\draw [fill=white] (-2.5,-5) circle (0.25);
\draw [draw=none, fill=lightgray] (-2.8,-4.85) rectangle (-3.35,-5.15);
\draw (-2.98,-4.85) -- (-3.25,-4.85);
\draw (-2.98,-5.15) -- (-3.25,-5.15);

\draw [fill=lightgray] (0,-7.5) circle (1);
\draw [fill=white] (0,-7.5) circle (0.75);
\draw [fill=lightgray] (0,-7.5) circle (0.5);
\draw [fill=white] (0,-7.5) circle (0.25);
\draw [draw=none, fill=lightgray] (-0.15,-7.8) rectangle (0.15,-8.35);
\draw (-0.15,-8.24) -- (-0.15,-7.97);
\draw (0.15,-8.24) -- (0.15,-7.97);

\draw [fill=lightgray] (0,-5) circle (1.25);
\draw [draw=none, fill=lightgray] (-1.2,-5.15) rectangle (-1.7,-4.85);
\draw (-1.22,-5.15) -- (-1.53,-5.15);
\draw (-1.22,-4.85) -- (-1.53,-4.85);
\draw [draw=none, fill=lightgray] (-0.15,-3.8) rectangle (0.15,-3.3);
\draw (-0.15,-3.47) -- (-0.15,-3.78);
\draw (0.15,-3.47) -- (0.15,-3.78);
\draw [draw=none, fill=lightgray] (1.2,-5.15) rectangle (1.7,-4.85);
\draw (1.22,-5.15) -- (1.53,-5.15);
\draw (1.22,-4.85) -- (1.53,-4.85);
\draw [draw=none, fill=lightgray] (-0.15,-6.2) rectangle (0.15,-6.7);
\draw (-0.15,-6.22) -- (-0.15,-6.53);
\draw (0.15,-6.22) -- (0.15,-6.53);
\draw (0,-5) node {$W$};
\draw (3.5,-7) node {$B\#L_1\#L_3\#L_5\#L_7$};
\draw (3.1,-6.7) -- (0.84,-5.92);

\end{tikzpicture}
}
\newcommand{\plotf}[1][]{
\begin{tikzpicture}[declare function={f(\x)=(100000^(-\x)*(2-\x^2)+100000^(\x-0.5)*exp(0.5-\x))/(100000^(-\x)+100000^(\x-0.5));}]
\begin{axis}[
axis lines = left,width=7cm,
  height=5.5cm,
ymin=0, ymax=2, ytick={0,1,2}, xmin=0, xmax=1,
 xtick={0,0.5,1},
 xticklabels={$0$,$\frac{1}{2}$,$1$},
 ylabel style={rotate=-90},
 x label style={at={(axis description cs:0.5,-0.1)},anchor=north},
    xlabel=$r$,
    ylabel={$f_1$}
    ]
\addplot [domain=0:1] {f(x)};
\end{axis}
\end{tikzpicture}
}
\newcommand{\plotg}[1][]{
\begin{tikzpicture}[declare function={g(\x)=(100000^(-\x)*(\x^2)+100000^(\x-0.5))/(100000^(-\x)+100000^(\x-0.5));}]
\begin{axis}[
axis lines = left,width=7cm,
  height=5.5cm,
ymin=0, ymax=1.5, xmin=0, xmax=1,ylabel style={rotate=-90},
 xtick={0,0.5,1},
 xticklabels={$0$,$\frac{1}{2}$,$1$},
 ytick={0,1},  
 yticklabels={$0$,$-\frac{T(p)}{2\pi}$},
 x label style={at={(axis description cs:0.5,-0.1)},anchor=north},
 y label style={at={(axis description cs:-0.03,0.5)},anchor=north},
    xlabel=$r$,
    ylabel={$f_2$}
    ]
\addplot [domain=0:1] {g(x)};
\end{axis}
\end{tikzpicture}
} 
\newcommand{\plotc}[1][]{
\begin{tikzpicture}

\draw [fill=lightgray] (-1,0) circle (2 and 1.5);
\path [fill=white] (-1,0.1) -- (-1,-0.2) arc [start angle=270, end angle=311, x radius=1, y radius=0.8] (-1,-0.2) arc [start angle=270, end angle=229, x radius=1, y radius=0.8] -- (-1,0.1);
\path [fill=white] (-1,0) -- (-1,0.2) arc [start angle=90, end angle=49, x radius=1, y radius=0.8] (-1,0.2) arc [start angle=90, end angle=131, x radius=1, y radius=0.8] -- (-1,0);

\path [fill=lightgray] (1.5,0) -- (2.1,0.7) arc [start angle=260, end angle=230, radius=1.21];
\path [draw=white, fill=white] (1.8,0) arc [start angle=180, end angle=90, x radius=0.3, y radius=0.7] (2.1,0.7) arc [start angle=260, end angle=230, radius=1.21];

\path [fill=lightgray] (1.5,0) -- (2.1,-0.7) arc [start angle=100, end angle=130, radius=1.21];

\draw (-1,-0.2) arc [start angle=270, end angle=330, x radius=1, y radius=0.8];
\draw (-1,-0.2) arc [start angle=270, end angle=210, x radius=1, y radius=0.8];
\draw (-1,0.2) arc [start angle=90, end angle=49, x radius=1, y radius=0.8];
\draw (-1,0.2) arc [start angle=90, end angle=131, x radius=1, y radius=0.8];

\path [draw=none, fill=gray] (2.1,0.7) arc [start angle=260, end angle=250, radius=1.21] arc [start angle=100, end angle=260, x radius=0.38, y radius=0.77];
\path [draw=none, fill=gray] (2.1,-0.7) arc [start angle=100, end angle=110, radius=1.21] -- (0.8,0);

\draw [draw=none, fill=gray] (2.1,0) circle (0.3 and 0.7);

\draw [fill=gray] (4.7,0) circle (1.84 and 1.84);

\draw (2.1,0) circle (0.3 and 0.7);

\draw [draw=none, fill=gray] (2.19,0.2) rectangle (2.9,0.4);
\draw [draw=none, fill=gray] (2.3,-0.4) rectangle (3,-0.2);

\draw [draw=none, fill=lightgray] (1.995,-0.664) arc [start angle=-70, end angle=70, x radius=0.3, y radius=0.7];
\draw [draw=none, fill=lightgray] (2.005,0.664) arc [start angle=110, end angle=250, x radius=0.3, y radius=0.7];

\draw (2.1,0.7) arc [start angle=260, end angle=230, radius=1.21];
\draw (2.1,-0.7) arc [start angle=100, end angle=130, radius=1.21];

\draw (2.2,0) arc [start angle=0, end angle=70, x radius=0.3, y radius=0.7];
\draw (2.2,0) arc [start angle=0, end angle=-70, x radius=0.3, y radius=0.7];

\draw [draw=none, fill=lightgray] (-0.2,-0.95) rectangle (1.55,0.95);
\draw (0.52,-0.96) -- (1.55,-0.96);
\draw (0.52,0.96) -- (1.55,0.96);

\draw [fill=lightgray] (4.7,0) circle (1.5 and 1.5);

\draw [fill=white] (4.7,0) circle (0.4);

\draw [draw=none, fill=lightgray] (2.1,-0.2) rectangle (3.4,0.2);

\draw (2.343,0.4) -- (2.9,0.4);
\draw (2.35,-0.4) -- (2.9,-0.4);
\draw (2.19,0.2) -- (3.2,0.2);
\draw (2.19,-0.2) -- (3.2,-0.2);

\draw (0.72,-1.3) node {$T_1$};
\draw (-0.5,-0.8) node {$X_1$};
\draw (1.47,-1.3) node {$T_2$};
\draw (4,-0.5) node {$X_3$};
\draw (1.07,0) node {$X_2$};
\draw (2,1) node {$B$};

\path [fill=lightgray] (1.3,0) -- (1.8,-0.7) arc [start angle=106, end angle=112, radius=1.21];
\draw (0.72,-0.96) -- (0.72,0.96);
\draw (1.47,-0.96) -- (1.47,0.96);
\draw [draw=none, fill=lightgray] (1.5,-0.2) rectangle (1.6,0.2);

\draw (1.89,0.757) arc [start angle=101, end angle=259, x radius=0.38, y radius=0.77];

\end{tikzpicture}
}

\begin{document}

\title{{Reeb flows without simple global surfaces of section}}
\author{Juno Kim, Yonghwan Kim, Otto van Koert}

\maketitle

\begin{abstract}
We construct, for any given positive integer $n$, Reeb flows on contact integral homology 3-spheres which do not admit global surfaces of section with fewer than $n$ boundary components.
We use a connected sum operation for open books to construct such systems.
Moreover, we prove that this property is stable with respect to $C^{4+\epsilon}$-small perturbations of the Hamiltonian given on the symplectization. 
\end{abstract}

\section{Introduction}
Global surfaces of section are an important tool to reduce the dynamics of flows on 3-manifolds to the dynamics of surface diffeomorphisms.
According to Ghys \cite{Ghys-09} global surfaces of section in their simplest form are a paradise for dynamicists.
In fact, just the existence of global surfaces of section gives information about the flow; 
for example, having a global surface of section without boundary for a flow on a $3$-manifold $M$ implies that $M$ fibers over the circle.
Since many $3$-manifolds, such as for example $S^3$, do not fiber over the circle, any global surface of section for such a manifold must have boundary, and hence the flow will have periodic orbits.

On the other hand, in \cite{Kuperberg94} Kuperberg has shown that there are flows on $S^3$ without any periodic orbits.
It is hence natural to look for global surfaces of section in a more restricted class of vector fields.
One candidate is the class of Hamiltonian vector fields, but without further restrictions, there are still difficulties.
For example, the horocycle flow on the unit cotangent bundle $ST^*\Sigma_g$ of a higher genus surface provides an example of a Hamiltonian flow on a compact manifold without any periodic orbits, see \cite{Asselle}.
Since $ST^*\Sigma_g$ does not fiber over the circle, it does not admit a global surface of section.

In both of the above examples, the main obstruction is the existence of periodic orbits.
In the case of Reeb flows, this obstruction vanishes.
In \cite{Taubes} Taubes has proved the so-called Weinstein conjecture for contact $3$-manifolds, which asserts the existence of periodic Reeb orbits on compact contact manifolds. 
This takes care of this particular obstruction, and there are indeed various results known on the existence of global surfaces of section for contact $3$-manifolds.

One of the first results is due Hofer, Wysocki and Zehnder, \cite{HWZ}, who have shown that Reeb flows on a dynamically convex 3-sphere must have disk-like global surfaces of section. 
The condition of dynamical convexity means, roughly speaking, that the winding number of the linearized flow is sufficiently large.
Hryniewicz and Salom\~ao \cite{Hryniewicz-Salomao} have extended this result by weakening the condition of dynamical convexity.

The global surfaces of section in these results have the simplest possible topology, and this can be useful in the analysis of the return map. This topology does not need to be simple.
Indeed, in \cite{Fried}, Fried has constructed global surfaces of section for transitive Anosov flows, pointing out that his construction methods do not control the genus of the surface of section.
Recently, existence results for global surfaces of section have been generalized to Reeb flows of $C^\infty$-generic contact forms by Contreras and Mazzucchelli in \cite{CM} and also by Colin, Dehornoy, Rechtman and Hryniewicz in \cite{CDRH}. Like Fried's results, these recent results also do not control the topology of the global surface of section.

We may thus view the simplest possible topology of a global surface of section on a fixed $3$-manifold as a measure of dynamical complexity of the flow. We note here that this is only a very partial measure of complexity: in general, it cannot distinguish between integrable and non-integrable flows, for instance.

In this paper, we will construct Reeb flows for which any global surface of section must have a much more ``complicated'' topology than those in \cite{HWZ} or \cite{Hryniewicz-Salomao}.
This extends the result of \cite{vK19} in which it was shown that there are Reeb flows without a disk-like global surface of section.
The precise statement of our first result is as follows:
\begin{thm}
\label{thm1}
Let $M$ be an integral homology 3-sphere with a contact structure $\xi$. 
Then for any integer $n>1$, there exists a contact form $\alpha$ for $\xi$ whose Reeb flow does not admit global surfaces of section with fewer than $n$ boundary  components.
\end{thm}

Now one may also ask whether these Reeb flows with ``complicated'' global surfaces of section are, in some sense, ``generic''. We answer this question positively by showing that our construction of the Reeb flow is stable under $C^
{\infty}$-small perturbations.
To keep the statement simple, we remind the reader that the Reeb flow from Theorem~\ref{thm1} is a Hamiltonian dynamical system on $(\R_{>0} \times M,d(\rho \alpha))$ with Hamiltonian $H=\rho$.

\begin{thm}
\label{thm2}
Let $M$ be an integral homology $3$-sphere with a contact structure $\xi$, and let $\alpha$ denote the contact form obtained in the proof of Theorem~\ref{thm1}.
Then there is a deformation $\bar H$ of the Hamiltonian $H=\rho$ on the symplectic manifold $(\R_{>0} \times M,d(\rho \alpha))$ with the following property:
For any $C^{4+\epsilon}$-small perturbation $\bar H_\delta$ of $\bar H$, the dynamics on the level set $\bar H_\delta=1$ do not admit global surfaces of section with fewer than $n$ boundary components.
\end{thm}
Complementary to the negative results here, one can also investigate what topologies are possible for a given Reeb flow.
For progress on this question, we mention the work of Albach, Geiges \cite{AG}, and that of Albers, Geiges and Zehmisch \cite{AGZ}.

To construct the Reeb flow for Theorem~\ref{thm1}, we will use the open book decomposition, and the book-connected sum operation. 
Since a connected sum with the standard $3$-sphere does not change the contact structure on $M$, performing book-connected sums with annuli will complicate the dynamics of the Reeb flow, while leaving the contact structure on $M$ unchanged. 
This is discussed in more depth in Section~\ref{sec:general_theory}, with an analysis of the invariant sets for the Reeb flow after the book-connected sum. 
In Section~\ref{secpf1}, we generalize the result from  \cite{vK19} by constructing a Reeb flow originating from an open book decomposition of $M$ which has many periodic orbits that do not link with each other. 
Then a linking number argument will show that a global surface of section for the Reeb flow must have many boundary components.

Theorem~\ref{thm2} will be proved in Section~\ref{sec:stability}.
To see that we still have enough invariant sets to run a linking argument we apply KAM theory.
This complicates the original construction somewhat, and makes the initial perturbation necessary.
Since the perturbation is in particular $C^2$-small, we still reconstruct the open book decomposition, and hence a global surface of section, on the perturbed level set.
The key point, stability near the binding, is explained in Appendix~\ref{appendix:stab_ob}; this only requires the usual implicit function theorem.
The global surface of section simplifies the analysis of the dynamics.

\section{General Theory}
\label{sec:general_theory}
\subsection{Global Surfaces of Section}
\noindent
Let $M$ be an oriented $3$-manifold, and let $\phi_t$ be a flow on $M$ generated by a vector field $X$ on $M$.
\begin{defn}
	A $\textit{global surface of section}$ for $(M,\phi_t)$ is a connected, compact, oriented surface $S$ embedded in $M$ such that the following conditions hold;
	\begin{enumerate}
		\item $X$ is positively transverse to the interior of $S$.
		\item For every $p\in M$, there exists some $t^+>0$ and $t^-<0$ such that $\phi_{t^+}(p),\:\phi_{t^-}(p)\in S$.
		\item The boundary of $S$ consists of periodic orbits of $\phi_t$.
	\end{enumerate}
\end{defn}
As indicated in the introduction, an important dynamical significance of this concept is that it allows us to convert problems of flows in $3$-manifolds into problems of surface diffeomorphisms.

We also get an immediate topological consequence from the definition of a global surface of section.
Namely, each orbit of $\phi_t$ lies either in the boundary of $S$, or positively intersects $S$. This leads us to consider the linking number of the periodic orbits.

As a reminder to the reader, the linking number of two oriented knots $k, \ell$, or more generally oriented links, in an oriented integral homology 3-sphere $M$ is defined as the oriented count of intersections between $\ell$ and a Seifert surface $F_k$ for $k$. 
This number is well-defined, i.e.~it is independent of the choice of Seifert surface.
We explain this in Appendix~\ref{app:seifert}.
We will also need some additional facts.
First of all, linking numbers are symmetric.
\begin{lemma}
	For two knots $\ell_1, \ell_2$ in an (oriented) integral homology 3-sphere $M$, the linking number is symmetric, so $\lk(\ell_1, \ell_2)=\lk(\ell_2, \ell_1)$.
\end{lemma}
\begin{proof}
We first recall that any oriented 3-manifold is an oriented boundary of an oriented 4-manifold. 
Indeed, by the Lickorish-Wallace theorem, all oriented 3-manifolds can be obtained by surgery on a framed link $L$ on $S^3$. 
Now view $S^3$ as the boundary of $D^4$, and attach 2-handles with the framing given by the framed link $L$ to obtain the oriented $4$-manifold $N$ with oriented boundary $M$.

We can take a collar neighborhood $N_M$ of $M$ in $N$. 
Also, let $F_1, F_2$ be Seifert surfaces in $M$ that have boundaries $\ell_1, \ell_2$. 
We can push the interior of $F_2$ into $N_M$ to make $F_2'$ in $N_M$ such that $F_2'\cap M$ is $\ell_2$. Then the intersection of $F_1, F_2'$ will only happen in $M$, therefore $F_1\cap F_2'=F_1\cap \ell_2$.
Furthermore, the orientations match as well.
The surfaces $F_2$ and $F_2'$ are isotopic rel boundary, so the intersection number $[F_1] \cdot [F_2]$ is equal to the intersection number $[F_1] \cdot [F_2']$. We conclude that $\lk(\ell_1,\ell_2)=[F_1]\cdot [F_2]$.
Similarly, we can repeat this process with $F_1$, choosing a similarly isotoped surface $F_1'$ (rel boundary).
By symmetry of intersections in the 4-dimensional manifold $N$ we find
	\begin{equation*}
	\lk(\ell_1,\ell_2)=[F_1]\cdot [F_2']=[F_2]\cdot [F_1']= \lk(\ell_2,\ell_1).
	\end{equation*}
\end{proof}

We will often make use of the following observation when dealing with linking numbers of links.
\begin{thm}
	Let $F$ be a surface in $M$ with boundary components $L_1,\cdots,L_k$. 
Assume that $\gamma$ is a link that does not intersect with $L_1,\cdots,L_k$. If $\:\lk(\gamma,L_i)=0$ for all $i$, then $\gamma$ also has algebraic intersection number 0 with $F$. 
\end{thm}
\begin{proof}
	Recall that $H_1(M\setminus L_i)\cong \mathbb{Z}$, since $M$ is a homology sphere. Therefore we can compute the linking number of $\gamma$ with $L_i$ as the number $n$ that satisfies $[\gamma]=n\cdot \alpha$ for a preferred (oriented) generator $\alpha$ of $H_1(M\setminus L_i)$. Since $\lk(\gamma,L_i)=0$ for all $i$, we can see that $[\gamma]=0$ in $H_1(M\setminus L_i)$.
	
	Another way to define the linking number is through Poincar\'e duality: first take a Seifert surface $K_i$ for each link component $L_i$, whose existence will be reviewed in Appendix~\ref{app:seifert}. 
	We can take the fundamental class $[K_i]\in H_2(K_i,L_i)$, and map it with the natural map $H_2(K_i,L_i)\to H_2(M,L_i)\cong H^1(M\setminus L_i)$ to the Poincar\'e dual $D([K_i])$.
	We can define the linking number of $\gamma$ and $L_i$ to be $\langle D([K_i]),[\gamma] \rangle$, where $[\gamma]\in H_1(M\setminus L_i)$ is the homology class generated by $\gamma$. It can be shown that these two definitions are equivalent: see the book of Gompf and Stipsicz, \cite{Gompf}.
	
	Consider the intersection number $[F] \cdot [\gamma]$ in $H_*(M\setminus L)$. 
	The same argument shows that we can compute $[F] \cdot [\gamma]$ by $\langle D([F]),[\gamma]\rangle$ for $[\gamma]\in H_1(M\setminus L)$.
	Since we have shown that the homology class generated by $\gamma$ is zero in each $H_1(M\setminus L_i)$, the class $[\gamma]\in H_1(M\setminus L)$ is also zero. Therefore, we conclude $[F] \cdot [\gamma]=0$.
\end{proof}

Now let us consider a global surface of section $S$ for a flow $\phi_t$. 
Let $\gamma$ be a periodic orbit of $\phi_t$ which is not a cover of any of the boundary orbits of $S$. We can apply Theorem 2.3 to show that either $\gamma$ has intersection number $0$ with $S$, or some $\lk(\gamma,L_i)$ is nonzero. The definition of a global surface of section excludes the first case, so we conclude

\begin{lemma} \label{link}
	Let $S$ be a global surface of section for a Reeb flow $\phi_t$ on an integral homology sphere $M$, and denote by $L_1,\cdots,L_k$ the boundary components of $S$. If $\gamma$ is a periodic orbit of $\phi_t$ which is not any of the $L_i$, then at least one of the linking numbers $\lk(\gamma,L_i)$ is nonzero for some $i$ in $1,\ldots,k$.
\end{lemma}

\subsection{Open book decompositions}
We will make use of so-called open books to prove our theorems.
Here is the definition.
\begin{defn}
	An \textit{open book decomposition} for an oriented manifold $M$ is a pair $(B,\pi)$ that satisfies
	\begin{enumerate}
		\item $B$ is an oriented link in $M$.
		\item $\pi:M\setminus B\to S^1$ is a (smooth) fiber bundle such that the closure of each fiber $F_{\theta}=\pi^{-1}(\theta)$ has boundary $B$.
	\end{enumerate}
\end{defn}
We call $B$ the binding, and the closure of each $F_{\theta}$ the page of the open book.
We will sketch briefly how this purely topological concept is related to global surfaces of section.

Suppose that $S$ is a global surface of section for an oriented  $3$-manifold $M$ with a smooth flow $\phi_t$.
Define $B:=\partial S$. This is an oriented link.
For each $x\in M\setminus B$ we define the minimal forward return time as
$$
\tau_+(x):=\inf_{t>0, \phi_t(x)\in S} t
$$
and the minimal backward return time as
$$
\tau_-(x):=\inf_{t>0, \phi_{-t}(x)\in S} t.
$$
Since the flow is smooth, we see that $\tau_\pm$ are smooth functions on $M\setminus S$.
Define the map
\[
\pi: M\setminus B \longrightarrow S^1 =\R /\Z,\quad
x \longmapsto 
\begin{cases}
[\frac{\tau_-(x)}{\tau_+(x)+\tau_-(x)}] & \text{if }x \notin S,\\
[0] & x\in S. 
\end{cases}
\]
We note that this map is continuous, but it is not smooth in general.
This means that we have constructed a continuous open book for which the global surface of section $S$ is a single page.
The map $\pi$ can be smoothed, but we won't need this construction, so we will not go into the details.

There is an alternative description of an open book, which is more convenient for constructions.
This is the following.
\begin{defn}
	Define an $\textit{abstract open book}$ as a pair $(\Sigma,\phi)$ such that
	\begin{enumerate}
		\item $\Sigma$ is an oriented, compact surface with boundary.
		\item The \textit{monodromy} $\phi:\Sigma\to\Sigma$ is a diffeomorphism restricting to the identity near the boundary.
	\end{enumerate}
\end{defn}

There is a correspondence between open book decompositions and abstract open books.
To go from an abstract open book to an open book decomposition, we do the following.
\begin{itemize}
\item construct the mapping torus
$$
M(\Sigma,\phi):=\Sigma \times \R /(x,\theta) \sim (\phi(x),\theta-2\pi)
$$
\item put $B:=\partial \Sigma$, and set $M:=B\times D^2 \cup_\partial M(\Sigma,\phi)$
\item define 
\[
\pi: M\setminus B \longrightarrow S^1=\R /2\pi \Z,
\quad
x \longmapsto
\begin{cases}
[\theta] & \text{if }x=(b;r,\theta)\in B\times D^2,\\
[\theta] & \text{if }x=[(s,\theta)]\in M(\Sigma,\phi).
\end{cases}
\]
The map $\pi$ is a well-defined, smooth map.
\end{itemize}
As a result, we obtain a smooth manifold $M$ together with an open book decomposition.

Conversely, given an open book decomposition $(B,\pi)$ on $M$, we define $\Sigma$ to be a page of the open book $\Sigma := \overline{\pi^{-1}([0])}$.
The monodromy $\phi$ can be constructed by choosing a connection that is standard near the binding. This is explained in \cite[Chapter 4.4.2]{Geiges:contacttopology}.
In the setting of contact open books, one needs a symplectic connection to obtain the symplectic monodromy. 
Details of the latter monodromy construction can be found in \cite{vK17}.

\subsection{Hamiltonian and Reeb vector fields}
Suppose that $(M^{2n},\omega)$ is a symplectic manifold, so $\omega$ is a closed, non-degenerate $2$-form.
Then given any smooth function $H:M\to \R$, commonly referred to as Hamiltonian, we can define the \emph{Hamiltonian vector field} $X_H$ of $H$ by the formula
\begin{equation}
\label{eq:def_ham_vf}
\iota_{X_H}\omega = -dH.
\end{equation}
Hamiltonian flows preserve the symplectic form, and the Liouville measure, i.e.~
$$
\mathcal L_{X_H}\omega =0, \quad
\mathcal L_{X_H}\omega^n =0.
$$
In order to define the Reeb vector field, we need an additional concept.
A \emph{contact form} on $Y^{2n+1}$ is a $1$-form $\alpha$ such that $\alpha \wedge d\alpha^n \neq 0$.
Given a contact form, we define the \emph{Reeb vector field} of $\alpha$ as the vector field $R_\alpha$ satisfying
\begin{equation}
\label{eq:def_reeb_vf}
\iota_{R_\alpha} d\alpha=0, \quad \iota_{R_\alpha} \alpha =1.
\end{equation}
Reeb vector fields are special Hamiltonian vector fields. Namely, if $(Y^{2n+1},\alpha)$ is a manifold with contact form $\alpha$, then we can form the \emph{symplectization} $(\R_{>0} \times Y, d(\rho \alpha) \,)$, where $\rho$ is the coordinate on $\R$.
The symplectization is a symplectic manifold, and the Hamiltonian vector field of the Hamiltonian $\rho$ on the level set $\rho=1$ equals the Reeb vector field $R_\alpha$.

\subsection{Contact open books}
\label{sec:ob_form}
We now adapt this open book setup to the setting of contact manifolds and Reeb flows.
We will essentially copy the above construction for going from an abstract open book to an open book decomposition, but now in the presence of a geometric structure.
We will follow the Giroux construction.

First we strengthen the requirements to be able to get a contact structure.
\begin{itemize}
\item We require the compact, oriented surface $\Sigma$ with boundary to be a so-called Liouville domain.
This means that $\Sigma$ comes equipped with a Liouville form $\lambda$. 
This is a $1$-form $\lambda$ such that $d\lambda$ is an area-form, which induces the given orientation on $\Sigma$, and such that $\lambda$ induces the natural orientation\footnote{as defined by the outward pointing normal} on the boundary $\partial \Sigma$.
\item We require the monodromy $\phi$ to be an exact symplectomorphism that is the identity near the boundary of $\Sigma$.
This means that $\phi$ is an area-preserving diffeomorphism that satisfies $\phi^*\lambda=\lambda - dT$. 
\end{itemize}
We can and will assume that $T<0$.

We slightly modify the definition of the mapping torus. 
We will define $M(\Sigma, \phi)$ as the quotient of  $\Sigma\times\mathbb{R}$ by a relation $\sim$, where $(p,\theta)\sim(\phi(p),\theta+T(p))$.
This mapping torus is diffeomorphic to the one we defined earlier. 
The reason to deform the mapping torus like this comes from the following observation.
The contact form $\mu=d\theta+\lambda$ descends to our deformed mapping torus $M(\Sigma,\phi)$. 
The Reeb vector field for this contact form $\mu$ is given by $\partial_\theta$.
In particular, we see that periodic Reeb orbits in the mapping torus correspond to periodic points of the monodromy $\phi$.

As before, we define a neighborhood of the binding of the open book decomposition as $B(\Sigma)= \partial\Sigma\times D^2$. We need to describe the gluing between a neighborhood of the binding and the mapping torus. 
Given an annulus $R(\frac{1}{2},1]$ of inner and outer radii $\frac{1}{2},1$, including the outer circle of radius 1, take
\begin{align*}
\Phi:\partial\Sigma\times R(\frac{1}{2},1]&\to (-\frac{1}{2},0]\times\partial \Sigma \times S^1\subset M(\Sigma,\phi), \\
(p,r,\theta)&\mapsto (\frac{1}{2}-r,p,-\frac{\theta T(p)}{2\pi}).
\end{align*}
Because $\phi$ is the identity near the boundary of $\Sigma$, the function $T$ is locally constant near the boundary.
Hence we see that $\Phi$ is a well-defined map from $\partial\Sigma\times R(\frac{1}{2},1]$ in $B(\Sigma)$ to $(-\frac{1}{2},0]\times\partial \Sigma\times S^1$ in $M(\Sigma,\phi)$. 
Define 
$$
M:=B(\Sigma)\amalg M(\Sigma,\phi) / \Phi,
$$
and write $B=\partial\Sigma\times\{0\}$ in $B(\Sigma)$. 
As before it is possible to construct a map $\pi$ making this into an open book decomposition on $M$, but we won't use this.

\begin{figure}
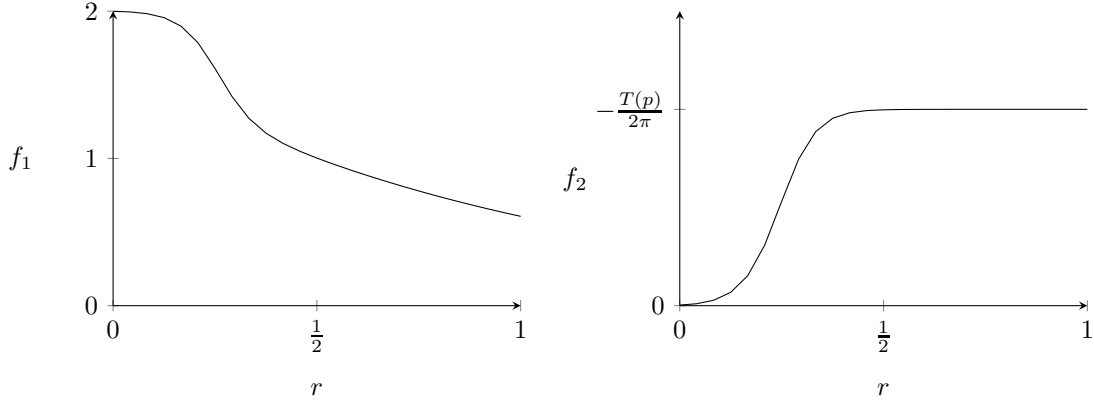

\centering
\plotf
\plotg
	\caption{The functions $f_1,f_2$ depending on $r$. \\$T(p)$ depends on the component of $B(\Sigma)$.}
	\label{plotfg}
\end{figure}

We already have a contact form $\mu=d\theta+\lambda$ on $M(\Sigma,\phi)$, which we will now extend to the whole space $M$.
Since $\lambda$ is a Liouville form, we can write $\mu$ near the boundary as $d\theta+e^{\tilde R}\lambda\vert_{\partial\Sigma}$, where ${\tilde R}$ is the collar parameter near the boundary.
This form is pulled back by $\Phi$ to $\Phi^*\mu=e^{\frac{1}{2}-r}\lambda\vert_{\partial\Sigma}-\frac{T(p)}{2\pi}d\theta$ on $B\times R(\frac{1}{2},1]$. 
To extend this form to all of $B(\Sigma)$, we will take smooth functions $f_1,f_2$ such that $\alpha=f_1\lambda\vert_{\partial\Sigma}+f_2d \theta$. 
We choose these functions $f_1,f_2$ to depend only on $r$, and to satisfy the following conditions:
\begin{enumerate}
	\item For $\frac{1}{2}\leq r \leq 1$, $f_1(r)=e^{\frac{1}{2}-r}$, $f_2(r) = -\frac{T(p)}{2\pi}$. 
	The latter is a constant depending on the component of $B(\Sigma)$. 
	\item For $r$ near 0, $f_1(r)=2-a r^2, f_2(r)=r^2$, where $a$ is a positive, irrational number to be determined later.
	\item $f_1f_2'-f_1f_2'>0$ for all $r>0$.
\end{enumerate}
The second condition ensures that $\alpha$ is a smooth $1$-form; the numerical constants are chosen to have a suitable Reeb flow near the binding for the construction in Section~\ref{sec:stability}.
Together with the third condition, this ensures that $\alpha$ is a positive contact form.
We sketch the graphs for $f_1,f_2$ in Figure~\ref{plotfg} above.

The upshot is that the manifold $M$ comes equipped with a contact form, which we will denote by $\alpha$.
We will denote the contact manifold constructed from the abstract open book $(\Sigma, \phi)$ by $\mathcal{OB}(\Sigma,\phi)$, and call it a contact open book.

\subsection{Reeb flow on a contact open book}
Near the binding of a contact open book $\mathcal{OB}(\Sigma,\phi)$ with contact form $\alpha$, the Reeb field has a particularly simple form.
Indeed, 
\begin{equation}
\label{eq:Reeb_vf}
R_\alpha=\frac{1}{f_1f_2'-f_2f_1'}(f_2'R_\lambda-f_1'\partial_\theta),
\end{equation}
where $R_\lambda$ is the Reeb flow for the contact form $\lambda|_{\Sigma}$.
In our case, where $\Sigma$ is just a surface, this vector field $R_\lambda$ is just the (positively oriented) unit tangent vector field to the boundary of $\Sigma$. 
The flow of this vector field has no component in the $\partial_r$ direction, so it preserves the $r$ coordinates in $\partial\Sigma\times D^2$. 
Hence we obtain the following invariant sets of this flow.
\begin{enumerate}
	\item The sets $\partial\Sigma\times S^1$ with fixed radius $r$ are invariant.
	\item The binding $B$ is an invariant set of the Reeb flow.
\end{enumerate}

From the explicit form of the flow, we see each page $F_\theta$ is a global surface of section for the Reeb flow.
Since every contact $3$-manifold admits a supporting open book by a result of Giroux, \cite{Giroux03}, it follows that given any contact $3$-manifold $(Y,\xi)$, there exists some Reeb flow that admits a global surface of section.
\begin{observation}
Every compact contact $3$-manifold $(Y,\xi)$ has some Reeb flow that admits a global surface of section.
\end{observation}
The stronger statement, whether every Reeb flow admits a global surface of section, can of course not be settled with the above techniques.
The following example of an open book decomposition will be a key ingredient for our constructions in the subsequent sections.

\begin{lemma} \label{s3}
	Let $A=[-1,1] \times S^1$ be an annulus equipped with the Liouville form $r d\theta$, and let $\phi:A\rightarrow A$ be a positive Dehn twist, defined as $\phi(x,\theta)=(x,\theta+\sigma(x))$ where $\sigma(x)=\pi(1-x)$. Then the manifold $\mathcal{OB}(A,\phi)$ carries the standard contact structure $\xi_0$ on $S^3$.
\end{lemma}
To see this:
This open book is actually the stabilization of the trivial open book $\mathcal{OB}(D^2,\mathrm{id})$ for $(S^3,\xi_0)$, and can be shown to be contactomorphic to $(S^3,\xi_0)$ using the methods from \cite{vK17}.
We review another operation on open books that will simplify our arguments later.

\begin{defn}
	Let $(\Sigma_1,\phi_1)$, $(\Sigma_2,\phi_2)$ be two abstract open books, and $c_0$, $c_1$ two arcs each properly embedded in $\Sigma_1$, $\Sigma_2$, and take rectangular neighborhoods $c_1\times [-1,1]$, $c_2\times[-1,1]$ each in $\Sigma_1$, $\Sigma_2$. The \textit{Murasugi sum} $\mathcal{OB}(\Sigma_1,\phi_1)*\mathcal{OB}(\Sigma_2,\phi_2)$ is defined as the manifold constructed from the abstract open book with page $\Sigma_1 \natural \Sigma_2$, where we identify $c_1\times\{-1,1\}$ with $\partial c_2\times[-1,1]$. The monodromy is defined as the composition $\tilde \phi_1\circ \tilde \phi_2$, where $\tilde \phi_1$ and $\tilde \phi_2$ are the extensions of $\phi_1$ and $\phi_2$ as the identity to the boundary connected sum $\Sigma_1 \natural \Sigma_2$. 
\end{defn}

A result of Torisu, \cite{Torisu} shows that $\mathcal{OB}(\Sigma_1,\phi_1)*\mathcal{OB}(\Sigma_2,\phi_2)$ is contactomorphic to $\mathcal{OB}(\Sigma_1, \phi_1)\#\mathcal{OB}(\Sigma_2,\phi_2)$. 
Note that the Murasugi sum reduces the number of boundary components of an open book. Thus, given any abstract open book for a manifold $M$, we may repeat this construction to obtain an abstract open book for $M$ with only one boundary component.

\subsection{Invariant sets for the book-connected sum construction}

We now review another operation on abstract open books, called the \textit{book-connected sum}. Let $(W_i,\psi_i)$, $i=1,2$ be two abstract open books, each with a contact structure. 
Then we can define a new abstract open book $(W_1\natural_{B_1,B_2}W_2,\psi_1\natural_{B_1,B_2}\psi_2)$ along specified boundary components $B_i\subseteq \partial W_i$. 
We will write $\natural$ for the boundary connected sum, an operation that can also be seen as $1$-handle attachment.
The symbol $\#$ stands for connected sum.
We use this notation both for operations on manifolds and for gluing maps together (silently extending a map as the identity if necessary).
The subscripts clarify where these operations are performed. We omit subscripts whenever their meaning is clear from the context. 

For completeness, here are our definitions.
The page $W_1\natural_{B_1,B_2}W_2$ is formed by attaching a Weinstein 1-handle $H$ to the disjoint union $W_1\coprod W_2$ along two Darboux balls in boundary, each in $B_1$ and $B_2$, respectively. 
The symplectomorphism $\psi_1\natural \psi_2$ is given by $\psi_i$ on the copy of $W_i$ in $W_1\coprod W_2$, and the identity on the handle. 
Note that we have the contact structure and Reeb vector field induced from the open book construction. Since the book-connected sum is a special case of the Murasugi sum we have explained earlier, we can apply the result of Torisu to show

\begin{lemma} \label{sums}
The book-connected sum $\mathcal{OB} (W_1\natural_{B_1,B_2}W_2,\psi_1\natural_{B_1,B_2}\psi_2)$ is contactomorphic to the contact connected sum $\mathcal{OB}(W_1,\psi_1)\#\mathcal{OB}(W_2,\psi_2)$.
\end{lemma}
Another argument for this lemma can be found in the appendix of \cite{vK17}.
The following description of the invariant sets will be useful.

\begin{lemma} \label{inv}
The induced Reeb flow decomposes the book-connected sum into four invariant sets: The handle orbits, the neighborhood $N(B_1\#B_2):=(B_1\#B_2)\times D^2$ of the boundary component used for the book-sum, and the remaining disjoint union of two ``page sets" $P_i$ described in the proof.
\end{lemma}

See Figure~\ref{booksum} for an example of such a decomposition.

\begin{proof}
The Reeb flow preserves the $r$-coordinate of the disk in the solid tori $\partial W\times D^2$, so $N(B_1\#B_2)$ is actually foliated by invariant tori. Since the monodromy is the identity on the handle, the orbits of $h\in H$ are the circles $\{h\}\times S^1$. Also, the mapping tori $M(W_i,\psi_i)$ viewed as subspaces of $M(W_1\natural W_2,\psi_1\natural \psi_2)$, are clearly disjoint invariant sets as the flows generated by each $\psi_i$ cannot cross into each other. Finally, the remaining sets are the neighborhoods $N(B)$ of the boundary components $B$ that are not the connect-summed knot $B_1\#B_2$, and thus are each foliated into invariant sets under the flow as well. Defining $P_i$ as the union 
$$
P_i:=M(W_i,\psi_i) \cup \bigcup_{B\text{ component of }\partial W_i \setminus B_i}N(B),
$$
the classification is complete.
\end{proof}
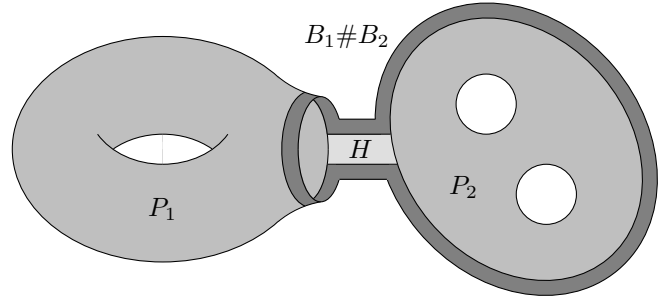
\begin{figure}[h]
\centering
\begin{tikzpicture}

\draw [fill=lightgray] (0,0) circle (2 and 1.5);
\path [fill=white] (0,0.1) -- (0,-0.2) arc [start angle=270, end angle=311, x radius=1, y radius=0.8] (0,-0.2) arc [start angle=270, end angle=229, x radius=1, y radius=0.8] -- (0,0.1);
\path [fill=white] (0,0) -- (0,0.2) arc [start angle=90, end angle=49, x radius=1, y radius=0.8] (0,0.2) arc [start angle=90, end angle=131, x radius=1, y radius=0.8] -- (0,0);

\path [fill=lightgray] (1.8,0) -- (2.1,0.7) arc [start angle=260, end angle=230, radius=1.21];
\path [draw=white, fill=white] (1.8,0) arc [start angle=180, end angle=90, x radius=0.3, y radius=0.7] (2.1,0.7) arc [start angle=260, end angle=230, radius=1.21];

\path [fill=lightgray] (1.8,0) -- (2.1,-0.7) arc [start angle=100, end angle=130, radius=1.21];
\path [draw=white, fill=white] (1.8,0) arc [start angle=180, end angle=270, x radius=0.3, y radius=0.7] (2.1,-0.7) arc [start angle=100, end angle=130, radius=1.21];

\draw (0,-0.2) arc [start angle=270, end angle=330, x radius=1, y radius=0.8];
\draw (0,-0.2) arc [start angle=270, end angle=210, x radius=1, y radius=0.8];
\draw (0,0.2) arc [start angle=90, end angle=49, x radius=1, y radius=0.8];
\draw (0,0.2) arc [start angle=90, end angle=131, x radius=1, y radius=0.8];

\path [draw=none, fill=gray] (2.1,0.7) arc [start angle=260, end angle=250, radius=1.21] arc [start angle=100, end angle=260, x radius=0.38, y radius=0.77];
\path [draw=none, fill=gray] (2.1,-0.7) arc [start angle=100, end angle=110, radius=1.21] -- (1.8,0);
\draw (1.89,0.757) arc [start angle=101, end angle=259, x radius=0.38, y radius=0.77];

\draw [draw=none, fill=gray] (2.1,0) circle (0.3 and 0.7);

\draw [fill=gray, rotate around={310:(4.7,0)}] (4.7,0) circle (2.1 and 1.7);

\draw [draw=none, fill=lightgray!50] (2.1,-0.2) rectangle (3.2,0.2);

\draw [draw=none, fill=lightgray] (1.995,-0.664) arc [start angle=-70, end angle=70, x radius=0.3, y radius=0.7];
\draw [draw=none, fill=lightgray] (2.005,0.664) arc [start angle=110, end angle=250, x radius=0.3, y radius=0.7];

\draw (1.8,0) arc [start angle=180, end angle=30, x radius=0.3, y radius=0.7];
\draw (1.8,0) arc [start angle=180, end angle=330, x radius=0.3, y radius=0.7];

\draw [draw=none, fill=gray] (2.19,0.2) rectangle (2.9,0.4);
\draw [draw=none, fill=gray] (2.3,-0.4) rectangle (3,-0.2);

\draw (2.1,0.7) arc [start angle=260, end angle=230, radius=1.21];
\draw (2.1,-0.7) arc [start angle=100, end angle=130, radius=1.21];

\draw (2.2,0) arc [start angle=0, end angle=70, x radius=0.3, y radius=0.7];
\draw (2.2,0) arc [start angle=0, end angle=-70, x radius=0.3, y radius=0.7];

\draw (2.343,0.4) -- (2.83,0.4);
\draw (2.35,-0.4) -- (2.98,-0.4);
\draw (2.19,0.2) -- (3.03,0.2);
\draw (2.19,-0.2) -- (3.12,-0.2);

\draw [fill=lightgray, rotate around={310:(4.7,0)}] (4.7,0) circle (1.9 and 1.5);

\draw [fill=white] (4.3,0.6) circle (0.4);
\draw [fill=white] (5.1,-0.6) circle (0.4);

\draw (0,-0.8) node {$P_1$};
\draw (4,-0.5) node {$P_2$};
\draw (2.63,0) node {$H$};
\draw (2.45,1.5) node {$B_1\#B_2$};

\end{tikzpicture}
\caption{Invariant sets of a book-connected sum: the page sets $P_i$, boundary neighborhood $N(B_1\#B_2)$ (dark grey), and the handle $H$.}
\label{booksum}
\end{figure}
A similar result clearly holds for multiple book-connected sums: each handle, boundary neighborhood, and page set are invariant sets, where the page set $P_i$ is now defined as the union of $M(W_i,\psi_i)$ and $N(B)$ for all boundary components $B$ of $\partial W_i$ not modified in the connected sum operations. This description of the invariant sets imply the following corollary:

\begin{corollary} \label{cor:orbit}
Let $p_1$, $p_2$ be periodic orbits in $P_1$, $P_2$, and let $h$ be a periodic orbit in the handle set $H$. Then $\lk(p_1,p_2)=\lk(p_1,h)=\lk(p_2,h)=0$.
\end{corollary}

\begin{proof}
We argue by finding a Seifert surface for $p_1$. By Lemma~\ref{sums}, $\mathcal{OB} (W_1\natural W_2,\psi_1\natural \psi_2)$ can be seen as a connected sum of $\mathcal{OB} (W_i,\psi_i)$, and we can assume that the balls used in the sum were contained in the solid tori $N(B_i)$. Since $P_i$ is invariant under the connected sum operation, we may consider $p_1$ as an orbit in $\mathcal{OB} (W_1,\psi_1)=P_1\cup N(B_1)$, and choose a Seifert surface $S\subset \mathcal{OB} (W_1,\psi_1)$ for $p_1$. Now we isotope $S$ so that $S$ intersects $N(B_1)$ only in disks of the form $\{v\}\times D^2$ for finitely many $v\in \partial W_1$. We further modify $S$ such that none of the $v$ are in the ball used for the connected sum. Thus, we may conversely view $S$ as a surface in $\mathcal{OB} (W_1\natural W_2,\psi_1\natural \psi_2)$, so that $S\subset P_1\cup N(B_1\# B_2)$. Finally, by the classification in Lemma~\ref{inv}, the surface $S$ cannot intersect either $p_2$ or $h$, so the corresponding linking numbers must be zero.
\end{proof}

\section{Proof of Theorem~\ref{thm1}} 
\label{secpf1}

Let $(M,\xi)$ be an oriented integral homology sphere with contact structure, and fix $n\geq 1$. Our objective is to construct a Reeb flow on $M$ which does not admit a global surface of section with $n$ or fewer boundary components. By a theorem of Giroux, \cite{Giroux03}, we can find an abstract open book $(W,\psi)$ that supports $\xi$, and we further assume $W$ to have a unique boundary component $B$. Let $\mathcal{OB}(A_i,\phi_i)$, $i=1,\cdots,n+1$ be copies of the open book of Lemma \ref{s3}, with $U_i=\{1\}\times S^1$, $L_i=\{-1\}\times S^1\subset A_i$ the boundary components of the annulus $A_i$. Consider the following book-connected sum:
\begin{align*}
X&=\mathcal{OB}(W\natural_{B,L_1}A_1\natural_{U_1,L_2}\cdots\natural_{U_{n},L_{n+1}}A_{n+1},\psi\natural_{B,L_1}\phi_1\natural_{U_1,L_2}\cdots\natural_{U_{n},L_{n+1}}\phi_{n+1})\\
&\cong\mathcal{OB}(W,\psi)\#_{B,L_1}\mathcal{OB}(A_1,\phi_1)\#_{U_1,L_2}\mathcal{OB}(A_2,\phi_2)\#_{U_2,L_3}\cdots\#_{U_{n},L_{n+1}}\mathcal{OB}(A_{n+1},\phi_{n+1})\\
&=\vcentcolon\mathcal{OB}_0\#\mathcal{OB}_1\#\cdots\mathcal{OB}_{n+1}\;\;\; \mathrm{in}\;\mathrm{shorthand}.
\end{align*}

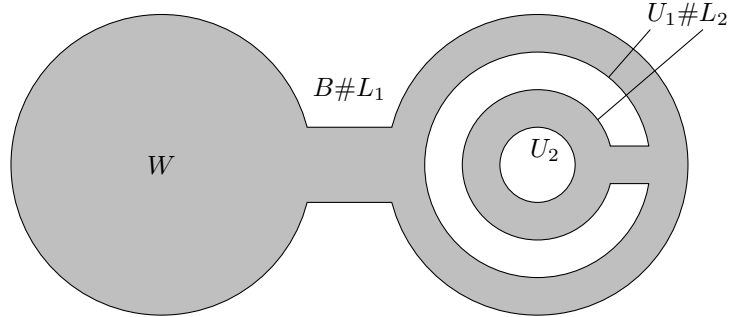
\begin{figure}[h!]
\centering
\begin{tikzpicture}
\draw [fill=lightgray] (-2.5,0) circle (2);
\draw [fill=lightgray] (2.5,0) circle (2);
\draw [draw=none, fill=lightgray] (-0.8,-0.5) rectangle (0.8,0.5);
\draw [fill=white] (2.5,0) circle (1.5);
\draw (-0.57,0.5) -- (0.57,0.5);
\draw (-0.57,-0.5) -- (0.57,-0.5);
\draw [fill=lightgray] (2.5,0) circle (1);
\draw [fill=white] (2.5,0) circle (0.5);

\draw [draw=none, fill=lightgray] (3.2,-0.25) rectangle (4.2,0.25);
\draw (3.46,-0.25) -- (3.985,-0.25);
\draw (3.46,0.25) -- (3.985,0.25);

\draw (-2.5,0) node {$W$};
\draw (2.6,0.2) node {$U_2$};
\draw (0,1) node {$B\#L_1$};
\draw (4.5,2) node {$U_1\#L_2$};
\draw (4,1.8) -- (3.45,1.17);
\draw (4.7,1.8) -- (3.3,0.6);

\end{tikzpicture}
\caption{A page of $X$ for $n=1$.}
\label{pf1}
\end{figure}

The page of $X$ for $n=1$ and $W=D^2$ is depicted in Figure~\ref{pf1}. Since the book-connected sum with the standard open book for $(S^3,\xi_0)$ does not change the contact structure, the open book $X$ is contactomorphic to the original homology sphere $(M,\xi)$.

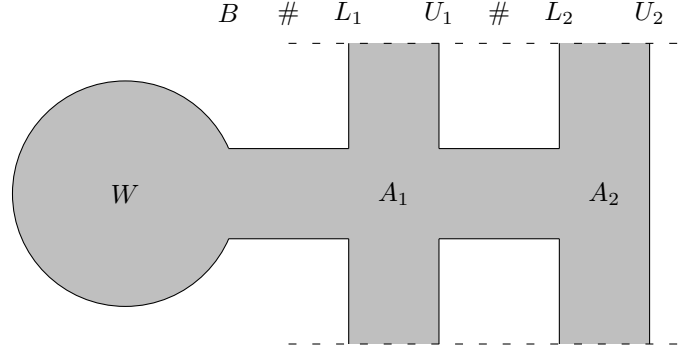
\begin{figure}
\centering
\begin{tikzpicture}
\draw [fill=lightgray] (-0.17,0) circle (1.5);
\draw [draw=none, fill=lightgray] (2.8,-2) rectangle (4,2);
\draw [draw=none, fill=lightgray] (5.6,-2) rectangle (6.8,2);
\draw [draw=none, fill=lightgray] (1.1,-0.6) rectangle (2.9,0.6);
\draw [draw=none, fill=lightgray] (3.9,-0.6) rectangle (5.7,0.6);
		
\draw (2.8,-2) -- (2.8,-0.6);
\draw (2.8,0.6) -- (2.8,2);
\draw (4,-2) -- (4,-0.6);
\draw (4,0.6) -- (4,2);
\draw (5.6,-2) -- (5.6,-0.6);
\draw (5.6,0.6) -- (5.6,2);
\draw (6.8,-2) -- (6.8,2);
\draw [loosely dashed] (2,2) -- (7.6,2);
\draw [loosely dashed] (2,-2) -- (7.6,-2);

\draw (1.2,-0.6) -- (2.8,-0.6);
\draw (1.2,0.6) -- (2.8,0.6);
\draw (4,-0.6) -- (5.6,-0.6);
\draw (4,0.6) -- (5.6,0.6);

\draw (1.2,2.4) node {$B$};
\draw (2,2.4) node {$\#$};
\draw (2.8,2.4) node {$L_1$};
\draw (4,2.4) node {$U_1$};
\draw (4.8,2.4) node {$\#$};
\draw (5.6,2.4) node {$L_2$};
\draw (6.8,2.4) node {$U_2$};
\draw (-0.17,0) node {$W$};
\draw (3.4,0) node {$A_1$};
\draw (6.2,0) node {$A_2$};
\end{tikzpicture}
\caption{An equivalent diagram for the page; the dashed lines are identified.}
\label{pf2}
\end{figure}

\begin{figure}[h]
\centering
\begin{tikzpicture}
\draw [draw=none, fill=lightgray] (1.33,-1.99) rectangle (10.395,1.99);

\draw [draw=none, fill=white] (1.2,-0.4) rectangle (2.8,0.4);
\draw [draw=none, fill=lightgray!50] (1.2,-0.4) rectangle (2.8,0.4);

\draw [draw=none, fill=gray] (-0.17,0) circle (1.5);

\draw [fill=lightgray] (-0.17,0) circle (1.3);

\path [fill=lightgray] (1.33,0) arc [start angle=0, end angle=15.5, radius=1.5] -- (0,0.4) -- (0,-0.1) -- (1.33,0) arc [start angle=0, end angle=-15.5, radius=1.5] -- (0,-0.4) -- (0,-0.1);

\draw (1.33,0) arc [start angle=0, end angle=20, radius=1.5];
\draw (1.33,0) arc [start angle=0, end angle=-20, radius=1.5];

\draw [draw=none, fill=gray] (1.2,0.4) rectangle (3,2);
\draw [draw=none, fill=gray] (1.2,-0.4) rectangle (3,-2);
\draw [draw=none, fill=gray] (3.8,0.4) rectangle (5.8,2);
\draw [draw=none, fill=gray] (3.8,-0.4) rectangle (5.8,-2);
\draw [draw=none, fill=gray] (6.6,0.4) rectangle (9.4,2);
\draw [draw=none, fill=gray] (6.6,-0.4) rectangle (9.4,-2);

\draw [draw=none, fill=white] (1.19,0.6) rectangle (2.8,2.1);
\draw [draw=none, fill=white] (1.19,-0.6) rectangle (2.8,-2.1);
\draw [draw=none, fill=white] (4,0.6) rectangle (5.6,2.1);
\draw [draw=none, fill=white] (4,-0.6) rectangle (5.6,-2.1);
\draw [draw=none, fill=white] (6.8,0.6) rectangle (9.2,2.1);
\draw [draw=none, fill=white] (6.8,-0.6) rectangle (9.2,-2.1);

\draw (-1.67,0) arc [start angle=180, end angle=23.5, radius=1.5];
\draw (-1.67,0) arc [start angle=180, end angle=336.5, radius=1.5];

\draw [draw=none, fill=white] (4,-0.4) rectangle (5.6,0.4);
\draw [draw=none, fill=lightgray!50] (4,-0.4) rectangle (5.6,0.4);
\draw [draw=none, fill=white] (6.8,-0.4) rectangle (9.2,0.4);
\draw [draw=none, fill=lightgray!50] (6.8,-0.4) rectangle (9.2,0.4);

\draw [draw=none, fill=white] (7.6,-1) rectangle (8.4,1);

\draw(1.06,-0.4) -- (3,-0.4);
\draw(1.06,0.4) -- (3,0.4);

\draw (2.8,-2) -- (2.8,-0.6);
\draw (2.8,0.6) -- (2.8,2);

\draw(3,-2) -- (3,-0.4);
\draw(3,0.4) -- (3,2);

\draw (4,-2) -- (4,-0.6);
\draw (4,0.6) -- (4,2);

\draw(3.8,-2) -- (3.8,-0.4);
\draw(3.8,0.4) -- (3.8,2);

\draw(3.8,-0.4) -- (5.8,-0.4);
\draw(3.8,0.4) -- (5.8,0.4);

\draw (5.6,-2) -- (5.6,-0.6);
\draw (5.6,0.6) -- (5.6,2);

\draw(5.8,-2) -- (5.8,-0.4);
\draw(5.8,0.4) -- (5.8,2);

\draw (6.8,-2) -- (6.8,-0.6);
\draw (6.8,0.6) -- (6.8,2);

\draw(6.6,-2) -- (6.6,-0.4);
\draw(6.6,0.4) -- (6.6,2);

\draw (1.2,-0.6) -- (2.8,-0.6);
\draw (1.2,0.6) -- (2.8,0.6);
\draw (4,-0.6) -- (5.6,-0.6);
\draw (4,0.6) -- (5.6,0.6);
\draw (6.8,-0.6) -- (7.6,-0.6);
\draw (6.8,0.6) -- (7.6,0.6);

\draw(6.6,-0.4) -- (7.6,-0.4);
\draw(6.6,0.4) -- (7.6,0.4);

\draw (8,0) node {$\mathbf{\cdots}$};

\draw (8.4,-0.6) -- (9.2,-0.6);
\draw (8.4,0.6) -- (9.2,0.6);

\draw(8.4,-0.4) -- (9.4,-0.4);
\draw(8.4,0.4) -- (9.4,0.4);

\draw (9.2,-2) -- (9.2,-0.6);
\draw (9.2,0.6) -- (9.2,2);

\draw(9.4,-2) -- (9.4,-0.4);
\draw(9.4,0.4) -- (9.4,2);

\draw (10.4,-2) -- (10.4,2);

\draw(2.8,-0.4) -- (2.8,0.4);
\draw(4,-0.4) -- (4,0.4);
\draw(5.6,-0.4) -- (5.6,0.4);
\draw(6.8,-0.4) -- (6.8,0.4);
\draw(9.2,-0.4) -- (9.2,0.4);

\draw [loosely dashed] (2,2) -- (11.2,2);
\draw [loosely dashed] (2,-2) -- (11.2,-2);

\draw (2,1.2) node {$N_0$};
\draw (4.8,1.2) node {$N_1$};
\draw (11,1.2) node {$U_{n+1}$};
\draw (-0.17,0) node {$P_0$};
\draw (3.4,0) node {$P_1$};
\draw (6.2,0) node {$P_2$};
\draw (9.8,0) node {$P_{n+1}$};
\draw (2,0) node {$H_0$};
\draw (4.8,0) node {$H_1$};

\end{tikzpicture}
\caption{The invariant sets of a page of $X$ for general $n$.}
\label{pf3}
\end{figure}
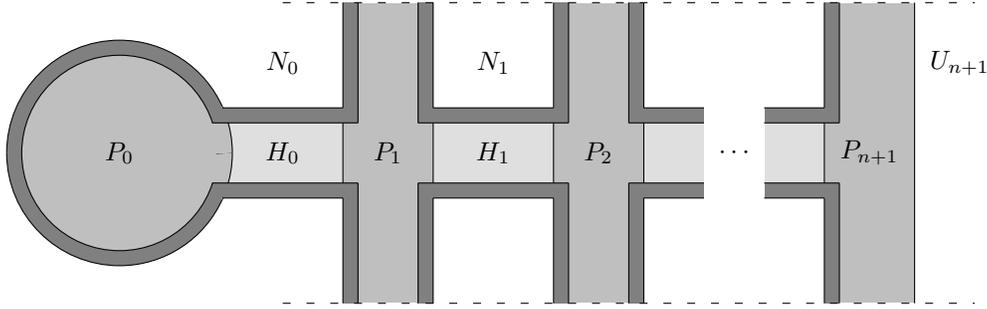

By viewing the annuli $A_i$ as rectangular strips with the ends glued together, we may draw an equivalent diagram which is convenient for the following argument: see Figure~\ref{pf2}. The 1-handles have been drawn at the same height for simplicity. We use the same description to illustrate the invariant sets in Figure~\ref{pf3}. Different colors each correspond to the page sets $P_i$, handle sets $H_i$, and boundary neighborhoods $N_i=N(U_i\#L_{i+1})$ and $N_0=N(B\#L_1)$.

The argument will proceed as follows: we will take a periodic orbit $h_i$ from each handle $H_i$. By Lemma~\ref{link}, each $h_i$ will either be the boundary of the global surface of section $S$, or will have positive linking number with at least one of the boundary components of $S$. However, Corollary~\ref{cor:orbit} strongly restricts the periodic orbits of nonzero linking number with $h_i$. We will use a counting argument to show that $\partial S$ should have more than $n$ components to satisfy the linking number condition. 
\par
Let us now go into the details. Assume that the induced flow on $X$ admits a global surface of section with $m(\leq n)$ boundary components $K_1,\cdots,K_m$. Since the handles consist of periodic orbits, we can find orbits $h_i$ in $H_i$ that are not equal to any of the $K_1,\cdots,K_m$. By Lemma~\ref{link}, the linking number $\lk(h_0,K_i)$ is nonzero for some $i$; without loss of generality, let $\lk(h_0,K_1)\neq 0$.
\par
In view of the sum $(\mathcal{OB}_0)\#(\mathcal{OB}_1\#\cdots)$, Corollary~\ref{cor:orbit} implies that only orbits in $H_0$ or $N_0$ can have nonzero linking number with $h_0$. In either case, $K_1$ is contained in the "left" page $P_0\cup H_0\cup N_0\cup P_1$ of the sum $(\mathcal{OB}_0\#\mathcal{OB}_1)\#(\mathcal{OB}_2\#\cdots)$. We apply  Corollary~\ref{cor:orbit} again to obtain $\lk(h_1,K_1)=0$. Therefore, $\lk(h_1,K_2)\neq 0$ for some $K_2$.
\par
We repeat this argument inductively. On the $i$th step, we have $\lk(h_j,K_{j+1})\neq 0$ for all $0\leq j\leq i$, so that $K_{j+1}\subset H_j$ or $N_j$. Therefore the knots $K_1,\cdots,K_{i+1}$ are all contained in the left-hand page of $(\mathcal{OB}_0\#\cdots\#\mathcal{OB}_{i+1})\#(\mathcal{OB}_{i+2}\#\cdots)$, so that none of these link with $h_{i+1}$. Accordingly, there must exist a boundary component $K_{i+2}$ such that $\lk(h_{i+1},K_{i+2})\neq 0$. However, at the $(m-1)$-st step, we cannot find a new boundary component $K_{m+1}$ with nonzero linking number with $h_m$. Therefore, we obtain the desired contradiction.
\qed


\section{Stability}
\label{sec:stability}

We now investigate if this boundary component condition for global surfaces of section is stable under perturbation, specifically of $C^\infty$-small deformations of the Hamiltonian on the symplectization $\R_{>0} \times \mathcal{OB}(W,\phi)$. In this section, we review Kolmogorov–Arnold–Moser (KAM) theory in order to study invariant sets of the perturbed Hamiltonian.
\par
We first state the Arnold-Liouville theorem. 
Recall that two functions $H_1, H_2$ on a symplectic manifold $(M,\omega)$ are said to be in involution if their Poisson bracket $\{H_1, H_2\}$ is zero. 
Given enough integrals in involution, the Arnold-Liouville theorem describes the invariant sets for the Hamiltonian motion.
\begin{thm}
	Let $(M^{2n},\omega)$ be a symplectic manifold of dimension $2n$. Let $H_1,\cdots,H_n$ be $n$ functions on $M$, which are in involution. 
	Suppose that $L$ is a compact, connected component of a regular level set of the function $H=(H_1,\cdots,H_n): M \to \R^n$. Then we have the following: 
\begin{itemize}
\item The set $L$ is a Lagrangian torus.
\item There is a neighborhood $\nu_M(L)$ that is diffeomorphic to $T^n\times D^n$ via the diffeomorphism
\[
\begin{split}
\Phi:T^n\times D^n \longrightarrow \nu_M(L).
\end{split}
\]
\item There are \emph{action-angle} coordinates $I, \phi$ on $D^n\times T^n$ with the properties
\begin{enumerate}
\item The symplectic form is standard:
$$
\Phi^*\omega=\sum_j dI_j\wedge d\theta^j.
$$
\item The coordinates $I_j$ depend only on the integrals $\{ H_k \}_{k=1}^n$.
\item The flow $\phi_{H_j}^t$ is linear in these coordinates, i.e. $\Phi^{-1}\circ \phi_{H_j}^t \circ \Phi(\theta,I)=(\theta+t \Omega_j(I),I)$.
\end{enumerate}
\end{itemize}	
\end{thm}
A proof can be found in Arnold \cite[Chapter 49]{Arnold} and Moser, Zehnder \cite[Section 3.1]{Moser_Zehnder}. We will call these invariant sets \textit{Liouville tori} for the given Hamiltonian. 
We will look at the Hamiltonian flow of the first Hamiltonian $H_1$ in these coordinates.
The above theorem tells us that the flow is linear on each Liouville torus, with slope $\Omega_1$.
We will just write $\Omega= \Omega_1$ and call it the frequency of a Liouville torus. 

The theory developed by Kolmogorov, Arnold and Moser, KAM, tells us that some of these Liouville tori survive perturbations.
We will review the statements from KAM theory that we will use and refer to \cite[Section 1,2]{Treschev} for proofs for the statements.

\begin{defn}
Let $(M^{2n},\omega)$ be a symplectic manifold equipped with a family of Hamiltonians $H=(H_1,\ldots, H_n)$. Assume that $T$ is a connected, compact, regular level set for which the Hamiltonians are in involution. Then $T$ is a Liouville torus associated with the Hamiltonian action of $H$. We will call this torus \textit{non-resonant} if $\Omega$ satisfies the condition that $k\cdot\Omega$ is nonzero for all $k\in\mathbb{Z}^n\setminus \{0\}$.

We call the first Hamiltonian $H_1$ \emph{non-degenerate} if the $n\times n$-matrix 
$$
(\frac{\partial^2 H_1}{\partial I_i \partial I_j})_{i,j}
$$
is invertible.
\end{defn}
In short, an integrable Hamiltonian $H_1$ is non-degenerate if the Jacobian of the frequency map is non-degenerate.
For non-degenerate (integrable) Hamiltonians, most Liouville tori are non-resonant tori.
\begin{thm}
Let $(M^{2n},\omega)$ be a symplectic manifold with a non-degenerate Hamiltonian $H=(H_1,\ldots,H_n)$. Then the non-resonant tori form a dense subset of the phase space.
\end{thm}

Now consider the case where we perturb the Hamiltonian. Given action-angle variables $I_1,\cdots,I_n$, $\theta^1,\cdots,\theta^n$ and a smooth function $F(I,\theta,\epsilon)$, we define a perturbation of the Hamiltonian to be $H(I,\theta, \epsilon)=H_0(I)+\epsilon F(I,\theta,\epsilon)$. To guarantee that these Liouville tori survive the perturbation, we will have to require a stronger condition of non-degeneracy. 

\begin{defn}
	A non-resonant Liouville torus with action-angle variables $I,\theta$ is said to be \textit{Diophantine} if the following inequality holds for some $c, \gamma>0$:
	\begin{equation*}
	\forall\: \kappa\in\mathbb{Z}^n\setminus\{0\}, \quad  \vert\langle \kappa, \Omega\rangle\vert\geq \frac{1}{c\lVert k\rVert^\gamma}.
	\end{equation*}	
\end{defn}

To indicate our strategy we first give a rough, imprecise version of the KAM theorem. We will state a more technical version later, Theorem~\ref{thm:twist_map}, that we actually use, see \cite[Theorem~2.1]{Treschev}
\begin{protothm}
\label{thm:KAM}
	Let $(M^{2n},\omega)$ be a symplectic manifold with non-degenerate Hamiltonian $H_0$. Let $I,\theta$ be the action-angle variables for $H_0$, and $f(I,\theta,\epsilon)$ a smooth function of sufficiently high regularity. 
	Define a small perturbation of the Hamiltonian $H(I,\theta,\epsilon)=H_0(I)+\epsilon f(I,\theta,\epsilon)$. 
	Then the following assertion holds:
	
\begin{itemize}
	\item A Diophantine torus in $M$ with respect to $H_0$ will survive a sufficiently small perturbation as a Diophantine torus with respect to $H$.
\end{itemize}
\end{protothm}

\subsection{Adaptation to our construction}
\label{sec:non-deg_Ham}
Let us outline how to apply this type of result to our setting of the symplectization $\R_{>0} \times \mathcal{OB}(\Sigma,\phi)$.
The main point of this section is to explain how we can ensure non-degeneracy of the Hamiltonian.
In general, we do not have an integrable system near the entire level set $\{ 1 \}\times M$, since the contact structure or topology on $M$ may obstruct the existence of integrals.
However, due to the special form of Reeb flow near the binding, namely the one from Equation~\eqref{eq:Reeb_vf}, we always have an integrable system near the binding, and all pieces of the book-connected sum also have such a structure.
With this in mind, let us start the construction.

Let $\alpha$ be the contact form constructed in the proof of Theorem~\ref{thm1}. Since the $1$-form $\lambda$ is just the standard angular form on the circle, we write
$$
\alpha=f_1(r)d\theta^1+f_2(r) d\theta^2.
$$
Let $H_0=\rho$ denote the Hamiltonian on $(\R_{>0} \times M, d(\rho \alpha) )$ which generates the Reeb flow.
Near the binding (but not at the binding), the action-angle coordinates are 
$$
I_1=\rho f_1(r), I_2 =\rho f_2(r), \theta^1, \theta^2.
$$
Indeed, we can compute the Jacobian of $I_1$ and $I_2$ to see that this is indeed a proper coordinate transformation for $\rho,r>0$.
We have
$$
\omega=d (\rho \alpha) =d I_1 \wedge d \theta^1 + d I_2 \wedge d \theta^2.
$$
With the special form of $f_1$ and $f_2$ chosen in Section~\ref{sec:ob_form}, we find that near the binding we have
$$
H_0 = \rho=\frac{1}{2}(I_1 +a I_2).
$$
This Hamiltonian is obviously degenerate with constant frequency $(\frac{1}{2},\frac{a}{2})$, but we can correct this by reparametrization and adding an integrable perturbation.
More explicitly we take the modified Hamiltonian
$$
\bar H_0= \frac{1}{4}(I_1 +a I_2)^2+g(r(I_1,I_2) ), \text{ where }
r(I_1,I_2)=\sqrt{
\frac{2I_2}{I_1+a I_2}
},
$$
where $g$ is a smooth function of $r$ that vanishes in a neighborhood of $r=0$.
The new frequency matrix is given by
$$
\Omega=
\left(
I_1+a I_2 +g' \frac{\partial r}{\partial I_1},
\quad
a(I_1+a I_2) +g' \frac{\partial r}{\partial I_2}
\right).
$$ 
Inspecting this expression, we can verify that we can make the new Hamiltonian $\bar H_0$ non-degenerate on the set $\rho=1$, and an open interval of $r$-values by choosing the function $g$ sufficiently general.
Furthermore, we can ensure that many Diophantine frequencies are attained. 
For example, we can go through the slope $\sqrt2$, which has continued fraction $1+\frac{1}{2+\frac{1}{2+\cdots}}$, and is Diophantine.
Below, we will apply Theorem~\ref{thm:KAM}.

\subsection{Stability of the global surface of section}
In this section, we will reconstruct a global surface of section for $C^{4+\epsilon}$-small perturbations of the Hamiltonian we have constructed above. This stability of the global surface of section will be used to prove stability of certain periodic orbits in the following section.
\par
Consider the contact manifold $(M,\alpha)$ defined by the book-connected sum $ X=\mathcal{OB}(W,\psi)\# \mathcal{OB}_1 \# \ldots \mathcal{OB}_{n+1}$ as in Section~\ref{secpf1}. The Reeb dynamics of the contact form $\alpha$ correspond to the Hamiltonian dynamics of $H_0=\rho$ on the symplectization $(\R_{>0} \times X, d (\rho \alpha)\, )$. In general, this is not an integrable system, since the monodromy of the open book on the homology sphere may obstruct the existence of first integrals.
We need to setup some notation to deal with this issue. 
We will write the space $X$ as $M \# M_0$, where $M_0 =\mathcal{OB}_1 \# \ldots \mathcal{OB}_{n+1}$ is contactomorphic to the tight contact $3$-sphere, where we keep in mind that a connected sum comes with the following decomposition
$$
M \# M_0 = M\setminus B \cup_\partial M_0 \setminus B_0,
$$
where $B$ is a Darboux ball in $M$ and $B_0$ is a Darboux ball in $M_0$.
Since the construction in the proof of Theorem~\ref{thm1} involves the book-connected sum, both $B$ and $B_0$ lie in a neighborhood of the binding, where first integrals do exist.
\par
To ``separate'' the non-integrable part of the dynamics, we define a cutoff function $\chi$ with the following properties:
\begin{itemize}
\item $\chi(x)=1$ for all $x\in X$ with $x\in M_0 \setminus B_0 \subset X$,
\item $\chi(x)=0$ for $x\in X$ with $x\in M(W,\psi) \subset M \setminus B \subset X$. Or in words, the function $\chi$ vanishes in the mapping torus region of $M$,
\item on the set where $\chi$ is not defined by the above, we can write $x=(\theta^1,r,\theta^2)\in \nu(B)=S^1\times D^2$ with angle coordinates $\theta^1$ and $\theta^2$.
We then choose $\chi$ to be a decreasing function of $r$ with the property that it equals $1$ for small $r$ (so as to be compatible with the first condition), and such that it vanishes for large $r$ (so that it is compatible with the second condition).
\end{itemize}
Now define 
$$
M_1:=\{x \in M ~|~\chi(x) =1 \}
.
$$
By restricting the Hamiltonian $H_0$ to a neighborhood of $\{ 1 \} \times M_1$ we obtain a completely integrable system with a complete flow.

\begin{figure}[H]
\centering
\begin{tikzpicture}

\draw [fill=lightgray] (0,0) circle (2 and 1.5);
\path [fill=white] (0,0.1) -- (0,-0.2) arc [start angle=270, end angle=311, x radius=1, y radius=0.8] (0,-0.2) arc [start angle=270, end angle=229, x radius=1, y radius=0.8] -- (0,0.1);
\path [fill=white] (0,0) -- (0,0.2) arc [start angle=90, end angle=49, x radius=1, y radius=0.8] (0,0.2) arc [start angle=90, end angle=131, x radius=1, y radius=0.8] -- (0,0);

\path [fill=lightgray] (2.8,0) -- (3.1,0.7) arc [start angle=265, end angle=244, radius=5.72];
\path [draw=white, fill=white] (2.8,0) arc [start angle=180, end angle=90, x radius=0.3, y radius=0.7] (3.1,0.7) arc [start angle=260, end angle=230, radius=1.21];

\path [fill=lightgray] (2.8,0) -- (3.1,-0.7) arc [start angle=95, end angle=116, radius=5.72];
\path [draw=white, fill=white] (2.8,0) arc [start angle=180, end angle=270, x radius=0.3, y radius=0.7] (3.1,-0.7) arc [start angle=100, end angle=130, radius=1.21];

\draw (0,-0.2) arc [start angle=270, end angle=330, x radius=1, y radius=0.8];
\draw (0,-0.2) arc [start angle=270, end angle=210, x radius=1, y radius=0.8];
\draw (0,0.2) arc [start angle=90, end angle=49, x radius=1, y radius=0.8];
\draw (0,0.2) arc [start angle=90, end angle=131, x radius=1, y radius=0.8];

\path [draw=none, fill=lightgray] (1.6,1) -- (3,0) -- (1.6,-1);

\path [draw=none, fill=gray] (3.1,0.7) arc [start angle=265, end angle=262.5, radius=5.72] arc [start angle=107, end angle=253, x radius=0.38, y radius=0.76];
\path [draw=none, fill=gray] (3.1,-0.7) arc [start angle=95, end angle=97.5, radius=5.72] -- (2.8,0);

\draw (2.85,0.73) arc [start angle=108, end angle=252, x radius=0.38, y radius=0.765];

\draw (2.3,0.83) arc [start angle=112, end angle=247, x radius=0.45, y radius=0.9];

\draw (1.7,1) arc [start angle=105, end angle=255, x radius=0.53, y radius=1.035];

\draw [draw=none, fill=gray] (3.1,0) circle (0.3 and 0.7);

\draw [draw=none, fill=lightgray] (2.995,-0.664) arc [start angle=-70, end angle=70, x radius=0.3, y radius=0.7];
\draw [draw=none, fill=lightgray] (3.005,0.664) arc [start angle=110, end angle=250, x radius=0.3, y radius=0.7];

\draw (2.8,0) arc [start angle=180, end angle=30, x radius=0.3, y radius=0.7];
\draw (2.8,0) arc [start angle=180, end angle=330, x radius=0.3, y radius=0.7];

\draw (3.1,0.7) arc [start angle=265, end angle=244, radius=5.72];
\draw (3.1,-0.7) arc [start angle=95, end angle=116, radius=5.72];

\draw (3.2,0) arc [start angle=0, end angle=70, x radius=0.3, y radius=0.7];
\draw (3.2,0) arc [start angle=0, end angle=-70, x radius=0.3, y radius=0.7];

\draw [fill=gray] (5.5,0) circle (1.7);
\draw [fill=lightgray] (5.5,0) circle (1.5);
\draw [fill=white] (5.5,0) circle (0.7);

\draw [draw=none, fill=gray] (3.2,0.2) rectangle (3.9,0.4);
\draw [draw=none, fill=gray] (3.2,-0.2) rectangle (3.9,-0.4);
\draw [draw=none, fill=lightgray] (3,-0.2) rectangle (4.5,0.2);

\draw (3.343,0.4) -- (3.85,0.4);
\draw (3.19,0.2) -- (4.01,0.2);
\draw (3.343,-0.4) -- (3.85,-0.4);
\draw (3.19,-0.2) -- (4.01,-0.2);

\draw (1.7,-1.35) node {$T_1$};
\draw (2.3,-1.2) node {$T_2$};
\draw (3,1) node {$B$};

\draw (0,-0.8) node {$X_1$};
\draw (1.65,0) node {$X_2$};
\draw [decorate,decoration={brace,mirror}] (2.3,-2) -- (7.2,-2) node[midway, below, yshift=-3pt]{$X_3$};

\end{tikzpicture}
\caption{Diophantine tori $T_1,T_2$ near the boundary component $B$ split the manifold into spaces $X_1,\:X_2,\:X_3$ with disjoint dynamics. The dark region represents the transit orbits between the homology sphere and the standard 3-sphere region. The global surface of section can be reconstructed for $C^{4+\epsilon}$-small perturbations of $\bar{H_0}$ in the $X_3$ region.}
\label{Separate}
\end{figure}
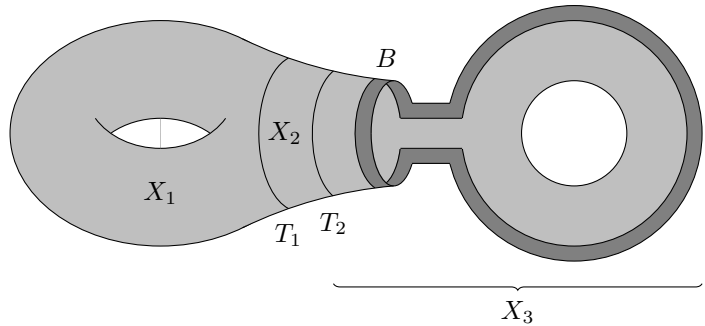

We now consider a $C^{4+\epsilon}$-small perturbation $\bar H_\delta:=\bar H_0 +\delta h$.
Since $\partial_\rho \bar H_0\neq 0$, we can apply the implicit function theorem and conclude that $\bar H_\delta^{-1}(1)$ is diffeomorphic to $\bar H_0^{-1}(1)$.
In fact, we see that $\bar H_\delta^{-1}(1)$ is a graph over $\{1\} \times M$, so we can still use $\chi$ to decompose this level set.
In particular we see that $\bar H_\delta$ defines a $C^{4+\epsilon}$-small perturbation of an integrable system near $\{1 \} \times M_1$.
Furthermore, we also see that $\bar H_\delta^{-1}(1)$ is contactomorphic to $\bar H_0^{-1}(1)$ by using Gray stability.
\par
We now apply Proposition~\ref{prop:stab_gss} from Appendix~\ref{appendix:stab_ob} to see that the level set $M_\delta:= \bar H_\delta^{-1}(1)$ admits a global surface of section $\Sigma_\delta$ for the Reeb flow on that level set. The situation is depicted in Figure~\ref{Separate}.
\subsection{Invariant orbits under small perturbations of $\bar H_0$}
We will now apply perturbation theory and KAM theory to identify two families of invariant sets. 
These invariant sets are the handle orbits corresponding to the critical points of the Morse Hamiltonian constructed in Appendix~\ref{app:morse_ham}, and the Diophantine tori in the annulus regions.

We first introduce a theorem from perturbation theory: the following statement is proved in the remarks from \cite[Theorem 2.2]{Moser_Zehnder}.

\begin{thm} 
\label{moser}
Consider an autonomous vector field $\dot{X}=f(X;\epsilon)$ on a smooth manifold $M$. 
Suppose that for $\epsilon=0$ there exists a periodic orbit $p(t;0)$, $p(0;0)=p$, of period $T>0$. Let $\varphi^T$ denote the time-$T$ flow of the system.

If 1 is a simple eigenvalue of $d\varphi_p^T:T_pM\rightarrow T_pM$, then for small $\epsilon$ there exists a periodic orbit $p(t,\epsilon)$ with period $T(\epsilon)$, such that $p(t,\epsilon)\rightarrow p(t,0)$ and $T(\epsilon)\rightarrow T$ as $\epsilon \rightarrow 0$. This orbit is unique up to a time shift.

In addition, suppose $\Sigma$ is a hypersurface transverse to the flow for small $\epsilon$, and let $\psi$ be the local diffeomorphism generated by the flow when $\epsilon=0$. Then the map $d\varphi_p^T$ has the following matrix representation:
\begin{gather*}
d\varphi_p^T=
\begin{pmatrix*}
1 & \cdots \\
0 & d\psi_p
\end{pmatrix*}.
\end{gather*}
\end{thm}

\begin{figure}[h]
	\centering
	\begin{tikzpicture}
\draw [draw=none, fill=lightgray] (0.6,-2) rectangle (6.05,2);
\draw [draw=none, fill=white] (1.5,1.5) arc [start angle=180, end angle=360, radius=0.5];
\draw [draw=none, fill=white] (1.5,-1.5) arc [start angle=180, end angle=0, radius=0.5];
\draw [draw=none, fill=white] (4,1.5) arc [start angle=180, end angle=360, radius=0.5];
\draw [draw=none, fill=white] (4,-1.5) arc [start angle=180, end angle=0, radius=0.5];
\draw [draw=none, fill=white] (1.5,2.01) rectangle (2.5,1.49);
\draw [draw=none, fill=white] (1.5,-2.01) rectangle (2.5,-1.49);
\draw [draw=none, fill=white] (4,2.01) rectangle (5,1.49);
\draw [draw=none, fill=white] (4,-2.01) rectangle (5,-1.49);
	
\draw (1.5,-2) -- (1.5,-1.5); \draw (1.5,2) -- (1.5,1.5);
\draw (1.5,1.5) arc [start angle=180, end angle=360, radius=0.5];
\draw (1.5,-1.5) arc [start angle=180, end angle=0, radius=0.5];
\draw (2.5,-2) -- (2.5,-1.5); \draw (2.5,2) -- (2.5,1.5);

\draw (4,-2) -- (4,-1.5); \draw (4,2) -- (4,1.5);
\draw (4,1.5) arc [start angle=180, end angle=360, radius=0.5];
\draw (4,-1.5) arc [start angle=180, end angle=0, radius=0.5];
\draw (5,-2) -- (5,-1.5); \draw (5,2) -- (5,1.5);

\draw (0.6,-2) -- (0.6,2); \draw (0.9,-2) -- (0.9,2); \draw (3.1,-2) -- (3.1,2); \draw (3.4,-2) -- (3.4,2);

\draw plot [smooth] coordinates {(1.2,-2) (1.3,-0.9) (2,0) (2.7,0.9) (2.8,2)};
\draw plot [smooth] coordinates {(1.2,2) (1.3,0.9) (2,0) (2.7,-0.9) (2.8,-2)};

\draw plot [smooth] coordinates {(3.7,-2) (3.8,-0.9) (4.5,0) (5.2,0.9) (5.3,2)};
\draw plot [smooth] coordinates {(3.7,2) (3.8,0.9) (4.5,0) (5.2,-0.9) (5.3,-2)};

\draw (5.6,2) -- (5.6,-2); \draw (5.9,2) -- (5.9,-2);

\draw (6.5,0) node {$\cdots$};

\draw [loosely dashed] (0.1,2) -- (6.5,2);
\draw [loosely dashed] (0.1,-2) -- (6.5,-2);
	\end{tikzpicture}
	\caption{Level sets of the Morse Hamiltonian. Note the critical points in each handle set, and the continuum of Liouville tori in each annulus set.}
	\label{level}
\end{figure}
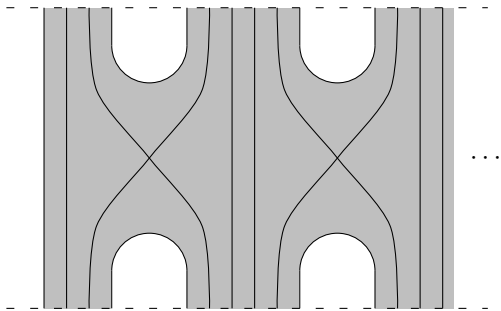

Now recall that the Morse Hamiltonian as constructed in Appendix~\ref{app:morse_ham} has a critical point at each handle set: the level sets are depicted in Figure~\ref{level}. We will apply the above theorem to show that these critical points correspond to invariant orbits with respect to the perturbation on $\bar H_0$.
\begin{cor} 
\label{handle}
After perturbation, there exists a periodic orbit $h_i$ on each handle connecting $A_i$ and $A_{i+1}$, near the original hyperbolic orbit, for $i=1,\ldots, 3n+1$.
\end{cor}

\begin{proof}
Take the hypersurface $\Sigma$ as a page of the open book before perturbation, and $p(t;0)$ as the original handle orbit. Let $p=p(0;0)=p(t;0)\cap \Sigma$. Since the return map $\psi$ is locally hyperbolic near the fixed point, $d\psi_p$ has no eigenvalue of modulus one. More precisely, near the center of the handle the Hamiltonian vector field constructed in Appendix~\ref{app:morse_ham} is up to rescaling $-2\pi (y \partial_x+x\partial_y)$, whose flow along $y=x$ (resp. $y=-x$) is expansion (resp. contraction) by $e^{2\pi t}$. The eigenvalues of $d\psi_p$ are therefore $e^{\pm 2\pi}$. Thus 1 is a simple eigenvalue of $d\varphi_p^T$ and the theorem applies.
\end{proof}

Clearly the orbits $h_i$ are unaffected by the handle attachment $Y=\mathcal{OB}_0 \# Y^0$, so we may view them as orbits in $Y$ with respect to the induced flow.
\par
In the proof of Theorem~\ref{thm1}, the handle sets functioned as `blockages' separating the page sets, so that page orbits cannot link with one another, and long orbits linking with multiple handle orbits cannot exist. 
Since we can now only ensure the existence of one orbit on each handle, this feature is lost. However, we are still able to obtain similar results by finding invariant tori which partition the manifold into regions with separate dynamics. The following invariant curve theorem is originally due to Moser 
\cite[Theorem 2.11]{Moser}, and its strengthening for lower regularity is due to Salamon, \cite{Salamon}.
\begin{thm}
\label{thm:twist_map}
Let $[a,b]\times S^1$ be an annulus with a twist mapping $(r,\theta)\mapsto (r,\theta+\sigma(r))$, such that $\sigma\in C^k$ for some $k>3$ and $|\sigma'|$ is bounded below by some positive constant. Then for any $\epsilon>0$, there exists $\delta>0$ such that all area-preserving mappings of the annulus into $\mathbb{R}^2$ of the form $$(r,\theta)\mapsto (f_1(r,\theta),\theta+f_2(r,\theta)),\quad \left\Vert f_1-r\right\Vert_{C^k}+\left\Vert f_2-\sigma \right\Vert_{C^k}< \delta$$ have an invariant curve of the following form, parametrized by $\gamma$:
$$ r=r_0+g_1(\gamma),\quad\theta=\gamma+g_2(\gamma)$$
where $g_1,g_2\in C^1$ and $\left\Vert g_1 \right\Vert_{C^1}+\left\Vert g_2 \right\Vert_{C^2}<\epsilon$. The induced mapping on the curve is given by $\gamma\mapsto \gamma+\kappa$, for some $\kappa$ incommensurable with $2\pi$.

Furthermore, for any choice of $\kappa \in\mathrm{im}\,\sigma$ satisfying the conditions
$$\left\vert \frac{\kappa}{2\pi}-\frac{p}{q}\right\vert \geq \alpha q^{-\beta},\quad \forall p,q\in\mathbb{Z},\;q>0$$
for some positive $\alpha,\beta$, there exists an invariant curve corresponding to $\kappa$ in the above sense.
\end{thm}

We now apply this theorem to the set $[-1/2,1/2]\times S^1$ in each annulus set $A_i$. Since the Hamiltonian vector field on this set is given by $\pi(1-r)\partial_\theta$, the twist condition is satisfied.
We perturb the Hamiltonian $\bar H_0$ by a $C^{4+\epsilon}$-small perturbation, so the Hamiltonian vector field is $C^{3+\epsilon}$-close to the unperturbed vector field, as is the return map. 
To guarantee that we have an annulus map, we apply a cutoff function $\delta$ to the perturbation.
This cutoff function vanishes near the boundary of the annulus and is $1$ in a smaller region in the annulus.
After that, the changes in the $r,\theta$ coordinates satisfy the condition for the Moser twist theorem. We can now choose $\kappa$ such that the invariant curve lies in the region where the cutoff function is $1$.

Therefore we can conclude that  
\begin{cor}
After perturbation, there is still an invariant curve $c_i$ on each annulus $A_i$, $i=1,\ldots,3n+2$, away from the upper and lower boundary circles.
\end{cor}
We will denote the invariant torus obtained by following the Reeb orbits through $c_i$ by $T_i$.

\section{Proof of Theorem 2}
We have identified invariant sets in both the handle and annulus sets in $X=\mathcal{OB}(W,\psi)\# \mathcal{OB}_1 \# \ldots \# \mathcal{OB}_{3n+2}$. 
We will first use some of the invariant tori to separate the dynamics of the homology sphere $M$ with the $M_0$ region. 
\par
\begin{figure}[h]
	\centering
	\begin{tikzpicture}
	\draw [draw=none, fill=lightgray] (1.3,-2) rectangle (10.4,2);
	
	\draw [draw=none, fill=gray] (-0.17,0) circle (1.5);
	
	\draw [fill=lightgray] (-0.17,0) circle (1.3);
	
	\path [fill=lightgray] (1.33,0) arc [start angle=0, end angle=15.5, radius=1.5] -- (0,0.4) -- (0,-0.1) -- (1.33,0) arc [start angle=0, end angle=-15.5, radius=1.5] -- (0,-0.4) -- (0,-0.1);
	
	\draw [draw=none, fill=gray] (1.2,0.4) rectangle (2,2);
	\draw [draw=none, fill=gray] (1.2,-0.4) rectangle (2,-2);
	\draw [draw=none, fill=gray] (2.8,0.4) rectangle (4.2,2);
	\draw [draw=none, fill=gray] (2.8,-0.4) rectangle (4.2,-2);
	\draw [draw=none, fill=gray] (5,0.4) rectangle (7.2,2);
	\draw [draw=none, fill=gray] (5,-0.4) rectangle (7.2,-2);
	\draw [draw=none, fill=gray] (8,0.4) rectangle (9.4,2);
	\draw [draw=none, fill=gray] (8,-0.4) rectangle (9.4,-2);
	
	\draw [draw=none, fill=white] (1.19,0.6) rectangle (1.8,2.1);
	\draw [draw=none, fill=white] (1.19,-0.6) rectangle (1.8,-2.1);
	\draw [draw=none, fill=white] (3,0.6) rectangle (4,2.1);
	\draw [draw=none, fill=white] (3,-0.6) rectangle (4,-2.1);
	\draw [draw=none, fill=white] (5.2,0.6) rectangle (7,2.1);
	\draw [draw=none, fill=white] (5.2,-0.6) rectangle (7,-2.1);
	\draw [draw=none, fill=white] (8.2,0.6) rectangle (9.2,2.1);
	\draw [draw=none, fill=white] (8.2,-0.6) rectangle (9.2,-2.1);
	
	\draw[draw=none, fill=lightgray] (1,-0.4) rectangle (2,0.4);
	
	\draw (-1.67,0) arc [start angle=180, end angle=23.5, radius=1.5];
	\draw (-1.67,0) arc [start angle=180, end angle=336.5, radius=1.5];
	
	\draw(1.06,-0.4) -- (2,-0.4);
	\draw(1.06,0.4) -- (2,0.4);
	
	\draw (1.8,-2) -- (1.8,-0.6);
	\draw (1.8,0.6) -- (1.8,2);
	
	\draw(2,-2) -- (2,-0.4);
	\draw(2,0.4) -- (2,2);
	
	\draw (3,-2) -- (3,-0.6);
	\draw (3,0.6) -- (3,2);
	
	\draw(2.8,-2) -- (2.8,-0.4);
	\draw(2.8,0.4) -- (2.8,2);
	
	\draw(2.8,-0.4) -- (4.2,-0.4);
	\draw(2.8,0.4) -- (4.2,0.4);
	
	\draw (4,-2) -- (4,-0.6);
	\draw (4,0.6) -- (4,2);
	
	\draw(4.2,-2) -- (4.2,-0.4);
	\draw(4.2,0.4) -- (4.2,2);
	
	\draw (5.2,-2) -- (5.2,-0.6);
	\draw (5.2,0.6) -- (5.2,2);
	
	\draw(5,-2) -- (5,-0.4);
	\draw(5,0.4) -- (5,2);
	
	\draw (1.2,-0.6) -- (1.8,-0.6);
	\draw (1.2,0.6) -- (1.8,0.6);
	\draw (3,-0.6) -- (4,-0.6);
	\draw (3,0.6) -- (4,0.6);
	\draw (4.2,-0.6) -- (4.2,-0.6);
	\draw (4.2,0.6) -- (4.2,0.6);
	
	\draw(5,-0.4) -- (7.2,-0.4);
	\draw(5,0.4) -- (7.2,0.4);
	
	\draw (5.2,0.6) -- (7,0.6);
	\draw (5.2,-0.6) -- (7,-0.6);
	
	\draw (7,0.6) -- (7,2); \draw (7,-0.6) -- (7,-2);
	\draw (7.2,0.4) -- (7.2,2); \draw (7.2,-0.4) -- (7.2,-2); 
	
	\draw (8.2,0.6) -- (8.2,2); \draw (8.2,-0.6) -- (8.2,-2);
	\draw (8,0.4) -- (8,2); \draw (8,-0.4) -- (8,-2); 
	
	\draw (8.2,-0.6) -- (9.2,-0.6);
	\draw (8.2,0.6) -- (9.2,0.6);
	
	\draw(8,-0.4) -- (9.4,-0.4);
	\draw(8,0.4) -- (9.4,0.4);
	
	\draw (9.2,-2) -- (9.2,-0.6);
	\draw (9.2,0.6) -- (9.2,2);
	
	\draw(9.4,-2) -- (9.4,-0.4);
	\draw(9.4,0.4) -- (9.4,2);
	
	\draw (10.4,-2) -- (10.4,2);
	
	\draw (3.2,0) node {$h_1$}; \draw [fill=black] (3.5,0) circle (0.04);
	\draw (5.4,0) node {$h_2$}; \draw [fill=black] (5.7,0) circle (0.04);
	\draw (9.3,0) node {$h_{3n+1}$}; \draw [fill=black] (8.7,0) circle (0.04);
	
	\draw (2.4,-2) -- (2.4,2); \draw (4.6,-2) -- (4.6,2); \draw (7.6,-2) -- (7.6,2);
	
	\draw (-0.17,-2.4) node {$\mathcal{OB}_0$};
	\draw (2.4,-2.4) node {$\mathcal{OB}_1$};
	\draw (4.6,-2.4) node {$\mathcal{OB}_2$};
	\draw (7.6,-2.4) node {$\mathcal{OB}_{3n+1}$};
	\draw (9.8,-2.4) node {$\mathcal{OB}_{3n+2}$};
	
	\draw (2.4,2.5) node {$\downarrow$};\draw (4.6,2.5) node {$\downarrow$}; \draw (7.6,2.5) node {$\downarrow$};
	\draw (2.4,2.9) node {$T_1$};\draw (4.6,2.9) node {$T_2$};\draw (7.6,2.9) node {$T_{3n+1}$};
	
	\path[draw,decorate,decoration={brace}] (2.45,2.2) -- (4.55,2.2);
	\path[draw,decorate,decoration=brace] (-1.6,2.2) -- (2.35,2.2);
	\path[draw,decorate,decoration=brace] (4.65,2.2) -- (7.55,2.2);
	\path[draw,decorate,decoration=brace] (7.65,2.2) -- (10.4,2.2);
	\draw (3.5,2.55) node {$V_1$}; \draw (0.375,2.55) node {$V_0$}; \draw (9.025,2.55) node {$V_{3n+1}$};
	
	\draw [draw=none, fill=white] (5.9,-1) rectangle (6.7,2.7);
	\draw (6.3,0) node {$\mathbf{\cdots}$};

	\draw [loosely dashed] (1,2) -- (11.2,2);
	\draw [loosely dashed] (1,-2) -- (11.2,-2);
	\end{tikzpicture}
	\caption{Invariant sets of a page of the perturbed flow.}
	\label{perturbpic}
\end{figure}
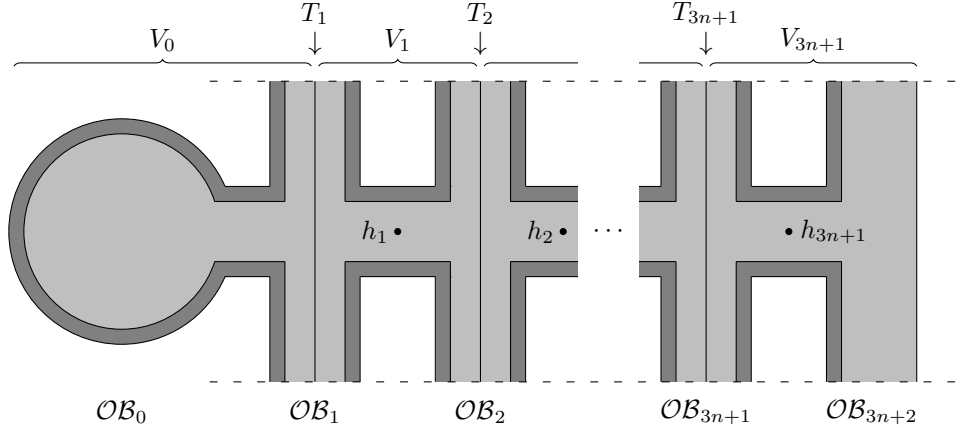

We label the invariant tori in each annulus $A_i$ as $T_i$, and each invariant handle orbit in the handle set $H_i$ to be $h_i$. Denote by $V_i$ the invariant set between $T_i$ and $T_{i+1}$: the situation is depicted in Figure~\ref{perturbpic}.

First, we will find two Diophantine tori $T_1,~T_1'$ in the tori ``connected'' to the homology sphere. As in Figure~\ref{Separate}, these Diophantine tori prevent the existence of orbits in the homology sphere that link with orbits in $\mathcal{OB}_2 \# \ldots \# \mathcal{OB}_{3n+2}$. Therefore, we can use a similar linking argument in the annuli connected sum region to provide a lower bound for the boundary components of the global surface of section.
\par

The key observation we make is the following:
\begin{proposition} \label{link2}
Orbits in $V_i$ and $V_j$ cannot link with each other if $|i-j|\geq 2$.
\end{proposition}

\begin{proof}
Consider the manifolds $Y^r=\mathcal{OB}_0\#\ldots\#\mathcal{OB}_{i+1}$ and $Y^{\ell}=\mathcal{OB}_{i+3}\#\cdots \# \mathcal{OB}_{3n+1}$. We may assume that the attachment of the $(i+1)$th handle was performed using a ball in $Y^r$ to the right of $T_{i+1}$ and a ball in $Y^{\ell}$ to the left of $T_{i+2}$. Thus, orbits in $V_i$ and $V_j$ are unaffected by the connected sum, and may be viewed as sitting in $Y^r$ and $Y^{\ell}$, respectively. Take a Seifert surface for an orbit $k$ contained in $V_i$. We can perform the book-connected sum along a Darboux ball that does not intersect the Seifert surface in $V_i$, which ensures the linking number with any other orbit contained in $V_j$ is zero.
\end{proof}

Now we can carry out the same linking number argument for the book-connected sum of $M$ with $3n+2$ copies of $S^3$ to conclude the proof. Assume that after a $C^{4+\epsilon}$ perturbation $h$, there exists a global surface of section on $\{1\}\times M$ with fewer than $n$ boundary components $K_1,\cdots,K_m$, for $m< n$. We first look at the handle orbits: each handle orbit $h_i$ should have positive linking number with a boundary component. Assume that $K_1$ has positive linking number with $h_2$. Then by Proposition~\ref{link2}, the knot $K_1$ cannot link with $h_l$ for $l\geq4$. We can repeat this process for each handle orbit to show that the global surface of section should have at least $n$ distinct boundary components, which yields a contradiction.

\begin{appendices}
\section{Seifert Surfaces for Integral Homology Spheres}
\label{app:seifert}
In this appendix, we will show that Seifert surfaces exist for any knot in an integral homology 3-sphere. Recall that a Seifert surface for an oriented link $k$ in an oriented 3-manifold $M$ is a connected oriented compact surface $S$ embedded in $M$ such that the oriented boundary $\partial S$ is equal to the link $k$.
\begin{thm}
Let $M$ be an integral homology 3-sphere, so $H_*(M;\Z)\cong H_*(S^3;\Z)$. 
For any oriented knot $k$ in $M$, there exists a Seifert surface $S$ for $k$ in $M$.
\end{thm}
\begin{proof}
Take a tubular neighborhood $N_k$ of $k$ in $M$, and let its boundary be $K$. Define $X$ to be the complement $M\setminus N_k$. 
We are going to construct the Seifert surface $S$ by first defining a map $f:K\to S^1$, and extending it to $\Tilde{f}:X\to S^1$. By the transversality theorem, we can assume that both $f, \Tilde{f}$ are transverse to some $p\in S^1$. 
Then $\Tilde{f}^{-1}(p)$ defines a surface $T$ with boundary in $K$. 
We then connect the boundary of $T$ in $K$ to $k$ to obtain a surface with boundary equal to $k$.
\par
First, we will define such a map $f$ from $K$ to $S^1$ whose preimage extends to a Seifert surface. For this reason, we need to identify $K$ with $S^1\times S^1$. 
There is a natural choice for a meridian $m$: it is a generator for $H_1(X)\cong \Z$. 
Then define a longitude for $K$ by taking a knot $l\in K$ such that intersection number in $K$ satisfies $ [l]\bullet [m]
=1$ and the homology class $[l]$ is trivial in the homology group $H_1(X,\Z)$. 
With these two knots, we can identify $K$ with $S^1\times S^1$. The knots $l$, $m$ give rise to global coordinates for $N_k$ as $(\theta_1,r,\theta_2)$ with $0\leq r\leq 1$ and $\theta_1,\theta_2\in S^1$, and its boundary $K$ will be the subset where $r$ takes the value $1$. We define $f:K\to S^1$ as the projection to the $\theta_2$ coordinate. 
\par
Now take the map $f$ defined above. 
Since $S^1$ is an Eilenberg-Maclane space $K(\Z,1)$, we can use a cohomology class $[f]$ in $H^1(K;\Z)$ to represent $f:K\to S^1$. Viewed as an equivalence class in De Rham cohomology, we can identify $[f]\in H_{DR}^1(K,\Z)$ to be $[d\theta_2]$.
We want to know whether the map $f$ extends to all of $X$.
For this, consider the inclusion $\iota:K\to X$. 
The extension of $f$ to a map $\tilde f:X \to S^1$ is equivalent to a cohomology class $[\Tilde{f}]$ in $H^1(X;\Z)$ that maps to $[f]\in H^1(K)$ under the induced map $\iota^*:H^1(X)\to H^1(K)$. From the cohomology exact sequence of the pair $(X,K)$ below, we can find such $[\Tilde{f}]$ if and only if $\partial^*[f]\in H^2(X,K)$ is zero:
\begin{equation*}
\begin{tikzcd}
H^1(X) \arrow[r,"\iota^*"] &H^1(K) \arrow[r,"\partial^*"] &H^2(X,K).
\end{tikzcd}
\end{equation*}
We will show that for every homology class $[S]\in H_2(X,K)$, the cohomology-homology pairing $\langle\partial^*[f],[S]\rangle$ takes the value $0$. For this purpose, we first identify the image of $\partial_*(H_2(X,K))$ in $H_1(K)$ in the relative long exact sequence.
Since $M$ is an integral homology sphere, we can see that the maps $\iota_1:K\xhookrightarrow{} N_k$, $\iota_2: K\xhookrightarrow{} X$ induce an isomorphism in the Mayer-Vietoris sequence:
\begin{equation*}
    (\iota_{1*},\iota_{2*}): H_1(K)\cong H_1(N_k)\oplus H_1(X).
\end{equation*}
Therefore, the image of the two generators $[l],[m]$ of $H_1(K)$ will generate $H_1(N_k)\oplus H_1(X)$. In particular, the image of the homology class $[m]$ in $H_1(X)$ will be a positive generator of $H_1(X)$, while the image of $[l]$ maps to zero in $H_1(X)$. Therefore, the image of $[l]$ in $H_1(N_k)$ will be a generator of $H_1(N_k)$. It follows that the kernel of the map $H_1(K)\to H_1(X)$ is generated by $[l]$. Now from the relative long exact sequence in homology
\begin{equation*}
\begin{tikzcd}
H_2(X,K) \arrow[r, "\partial_*"] &H_1(K) \arrow[r] &H_1(X),
\end{tikzcd}
\end{equation*}
the image of $\partial_*(H_2(X,K))$ in $H_1(K)$ is generated by $[l]$. In particular for any oriented piecewise smooth surface $S$ in $X$ with boundary in $K$, the homology class of the boundary $\partial S$ is an integer multiple of $[l]$. We look at the pairing of homology and cohomology to obtain the following equality
\begin{equation*}
    \langle\partial^*[f], [S]\rangle = \langle [f], \partial_* [S]\rangle =\langle[f],[\partial S]\rangle.
\end{equation*}
The homology class $[\partial S]$ in $H_1(K)$ is generated by $[l]$, which is represented by the coordinate $\theta_1$, while the cohomology class $[f]\in H^1(K)$ is identified with $[d\theta_2]$ in De Rham cohomology. Therefore we can conclude that for any oriented piecewise smooth surface $S$ in $X$ with boundary in $K$,
\begin{equation*}
    \langle\partial^*[f], [S]\rangle = 0,
\end{equation*}
which shows that $\partial^*[f]=0$. Therefore we can conclude that there is an extension $\Tilde{f}:X\to S^1$.

\par
From the transversality theorem, we can assume that a regular value $p\in S^1$ exists such that $f$, $\Tilde{f}$ are both transverse to $S^1$ at $p$. Therefore, the preimage $T$ defined as $f^{-1}(p)$ is a surface in $X$ with boundary in $K$. Now we look at the tubular neighborhood $N_k$ of $k$. The boundary of $T$ is in the boundary of $N_k$, so we can extend $T$ to a surface $S$ with boundary $k$ by connecting the boundary of $T$ to $k$ in $N_k$. Since $f^{-1}(p)$ is equal to the longitude $l$ in $H_1(K)$, the surface $S$ can be well-defined, possibly non-smooth at the boundary $K$ of $T$. By smoothing the surface at $K$, we can construct a Seifert surface $S$ for the knot $k$.
\end{proof}

Now recall that we have defined the intersection number of two oriented knots $k,l$ in an integral homology 3-sphere $M$ to be the intersection number between $l$ and the Seifert surface $F_k$ of $k$. We prove that this number is independent of the choice of Seifert surface for $k$.

\begin{proposition}
Assume that $F_1$, $F_2$ are two Seifert surfaces for an oriented knot $k$ in an integral homology sphere $M$. Then for any oriented knot $l$ in $M$, the intersection number of $l$ with $F_1$ agrees with the intersection number of $l$ and $F_2$.
\end{proposition}
\begin{proof}
A Seifert surface for $k$ is an oriented surface embedded in $M$ such that its boundary equals the oriented knot $k$. In terms of homology, the induced map $\partial_*:H_2(M,k)\to H_1(k)$ is an isomorphism because $M$ is a homology sphere.
Therefore the homology class $[F_1]-[F_2]$ is contained in the kernel of the map $\partial_*$, so $[F_1]=[F_2]$ in $H_2(M,k)$. It follows that the intersection number of $[l]$ with $F_1$ and $F_2$ agree.
\end{proof}

We have also implicitly used that Liouville tori inside integral homology 3-spheres divide the manifold into two connected components: we provide a short proof using the Mayer-Vietoris sequence.
\par
Assume a surface $T$ embedded in $M$ homeomorphic to the two-torus. Let $N$ be a tubular neighborhood of $T$ in $M$, with a homeomorphism $\phi:(-\epsilon,\epsilon)\times T\to N$. Then we can form two open sets $A=M\setminus T$, and $B=N$ that cover the total space $M$. 
We remark that $A\cap B$ can be identified with $\{(-\epsilon,0)\cup(0,\epsilon)\}\times T$ using the map $\phi$. 
The Mayer-Vietoris sequence for the pair $(A,\:B)$ gives a short exact sequence
\begin{equation*}
\begin{tikzcd}
H_1(M) \arrow[r] &H_0(T\times\{-\epsilon,\:\epsilon\}) \arrow[r] &H_0(A)\oplus H_0(B) \arrow[r] &H_0(M) \arrow[r] &0.
\end{tikzcd}
\end{equation*}
Since we can identify the homology groups in the sequence to be $H_0(T\times\{-\epsilon,\:\epsilon\})\cong\Z^2$, $H_0(B)\cong\Z$, $H_0(M)\cong\Z$, we can conclude that $H_0(A)$ has rank 2, and therefore that $M\setminus T$ has exactly 2 connected components.

\section{Construction of the Morse Hamiltonian}
\label{app:morse_ham}

In this appendix we will give an explicit construction of the Morse Hamiltonian $H$ on the page $Y^0$. The induced Hamiltonian flow will also have the following properties, which we use in the proof of Theorem~\ref{thm2}:
\begin{enumerate}
    \item $H$ has a critical point at each handle set, which corresponds to a hyperbolic periodic orbit for the induced Hamiltonian flow.
    \item The level sets of $H$ on the annulus part form a continuum of Liouville tori, some of which are Diophantine tori.
    \item The monodromy induced by the Hamiltonian flow is isotopic to the return map of $Y^0$, and is the identity near the boundary.
\end{enumerate}

We first present the construction for the union of an annulus and a handle set. Consider a standard model for such a set given as the domain $D=[-1,1]\times S^1\cup_\partial [1,\frac{5}{2}]\times[-\frac{\pi}{4},\frac{\pi}{4}]$, with the Liouville form $rd\theta$. The two sets $[-1,1]\times S^1$, $[\frac{1}{2},\frac{5}{2}]\times[-\frac{\pi}{4},\frac{\pi}{4}]$ each correspond to the annulus and handle sets. We assign to each domain a Hamiltonian $H_1, H_2$ given by
\begin{enumerate}
    \item $H_1$: $[-1,1]\times S^1\to \R: (r,\theta)\mapsto -\frac{\pi}{2}r^2+\pi r$,
    \item $H_2$: $[\frac{1}{2},\frac{5}{2}]\times[-\frac{\pi}{4},\frac{\pi}{4}]\to \R: (r,\theta)\mapsto C(1-(r-\frac{3}{2})^2+(\frac{4}{\pi}\theta)^2)$,
\end{enumerate}
where the constant $C>0$ will be determined later in the proof.

\begin{figure}
	\centering
	\begin{tikzpicture}
	\draw (-2,-2) rectangle (0,2);
	\draw [loosely dashed] (-3,-2) -- (1,-2);
	\draw [loosely dashed] (-3,2) -- (1,2);
	\fill[gray] (-2,-2) rectangle (0,2);
	\fill[white] (-1,-1) rectangle (0,1);
	\draw (0,-1) rectangle (1,1);
	\fill[lightgray] (0,-1) rectangle (1,1);
	\draw [white, line width=2pt] (0,-1) -- (0,1);
	\draw [line width=2pt, line cap=round] (0,-1) -- (-1,-1) -- (-1,1) -- (0,1);
	\draw [line cap=round] plot [smooth] coordinates {(0,0.9) (-0.8,0.8) (-0.9,0) (-0.8,-0.8) (0,-0.9)};
	\draw [line cap=round] plot [smooth] coordinates {(0,0.85) (-0.6,0.6) (-0.7,0) (-0.6,-0.6) (0,-0.85)};
	\draw [line cap=round] plot [smooth] coordinates {(0,0.8) (-0.4,0.5) (-0.5,0) (-0.4,-0.5) (0,-0.8)};
	\draw [line cap=round] plot [smooth] coordinates {(0,0.7) (-0.2,0.4) (-0.25,0) (-0.2,-0.4) (0,-0.7)};
	\draw [line width=2pt, line cap=round] (0,-0.65) -- (0,0.65);
	
	\draw (-2.7,0) node {$\rho=1$};
	\draw (1.7,0) node {$\rho=0$};
	\end{tikzpicture}
	\caption{Level sets of the cutoff function $\rho$ (not to scale). The left rectangle corresponds to the annulus set, and the right square corresponds to the handle set.}
	\label{cutoff}
\end{figure}
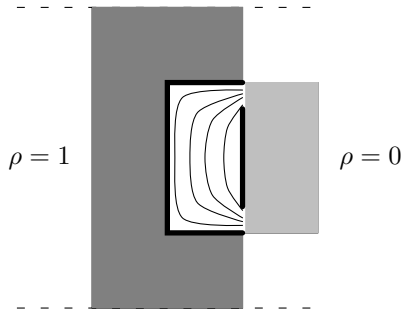

We now define a cutoff function to connect the level sets on the handle and the annulus. Take $\rho$ as a smooth function from $[\frac{1}{2},1]\times[-\frac{\pi}{4},\frac{\pi}{4}]$ to $[0,1]$ such that $\rho=1$ on $\{\frac{1}{2}\}\times [-\frac{\pi}{4},\frac{\pi}{4}]\cup [\frac{1}{2},1]\times \{\pm \frac{\pi}{4}\}$ and $\rho=0$ on $\{1\}\times [-\frac{\pi}{4},\frac{\pi}{4}]$. We extend $\rho$ to the whole set $D$ by assigning constant values $1$, $0$ such that $\rho$ is smooth. The level sets of $\rho$ are sketched in Figure \ref{cutoff}. Note that the cutoff function $\rho$ has discontinuities at the region where the annulus and handle sets attach.
\par
Now define the Hamiltonian $H_0=\rho H_1+(1-\rho)H_2$. The level sets of $H_0$ and the flow with respect to the Hamiltonian vector field is depicted in Figure~\ref{morse}. The constant $C$ is chosen such that the level sets match as in Figure~\ref{morse}. Note that this construction is made to trim the boundary to a smooth submanifold. Therefore we must check that the level sets of this Hamiltonian behave as in the Figure~\ref{morse}. We will use Morse theory arguments to determine the topology of the level sets. Since the only critical point of $H_0$ is contained in the handle set, the homotopy type of $H_0^{-1}(x)$ only changes when $x=-C$. Since the level set $H_0^{-1}(-C)$ behaves as in Figure~\ref{morse}, for small $\epsilon>0$,  $H_0^{-1}(-C+\epsilon)$ can be used to ``trim off'' the boundary to a smooth set. If we choose $\epsilon$ small enough, we can also make the level set $H_0^{-1}(-C+\epsilon)$ to not contain any discontinuities of $\rho$. Therefore, we can restrict $H_0$ to $H_0^{-1}(-C+\epsilon)$ as a smooth function.
\par
We now check if the constructed Hamiltonian satisfies our claimed conditions. The Hamiltonian vector field for $H_1, H_2$ can be computed to be $X_{H_1}=\pi(-r+1)\partial_\theta$, $X_{H_2}=-2C(\frac{4}{\pi}\theta\partial_r+(r-\frac{3}{2})\partial_\theta)$. Since the Hamiltonian vector field on the annulus part generates a positive Dehn twist, we can ensure that the contact manifold generated by the return map is $S^3$ with its standard tight contact structure. Therefore, the return map is isotopic to the return map of the book-connected sum.
\par
To check conditions (1),~(2), we will look at the level sets of $H_0$. On the handle region, the Hamiltonian $H_2$ has a hyperbolic critical point for $(r,\theta)=(\frac{3}{2},0)$, which corresponds to the hyperbolic periodic orbit. In the annulus region, the Hamiltonian vector field $X_1=\pi(-r+1)\partial_\theta$ generates Liouville tori for $-\frac{1}{2}\leq r\leq \frac{1}{2}$. Therefore, we have checked that conditions (1), (2) are satisfied.
\par

\tikzset{->-/.style={decoration={
			markings,
			mark=at position \halfway with {\arrow{Latex[length=2mm]}}},postaction={decorate}}}

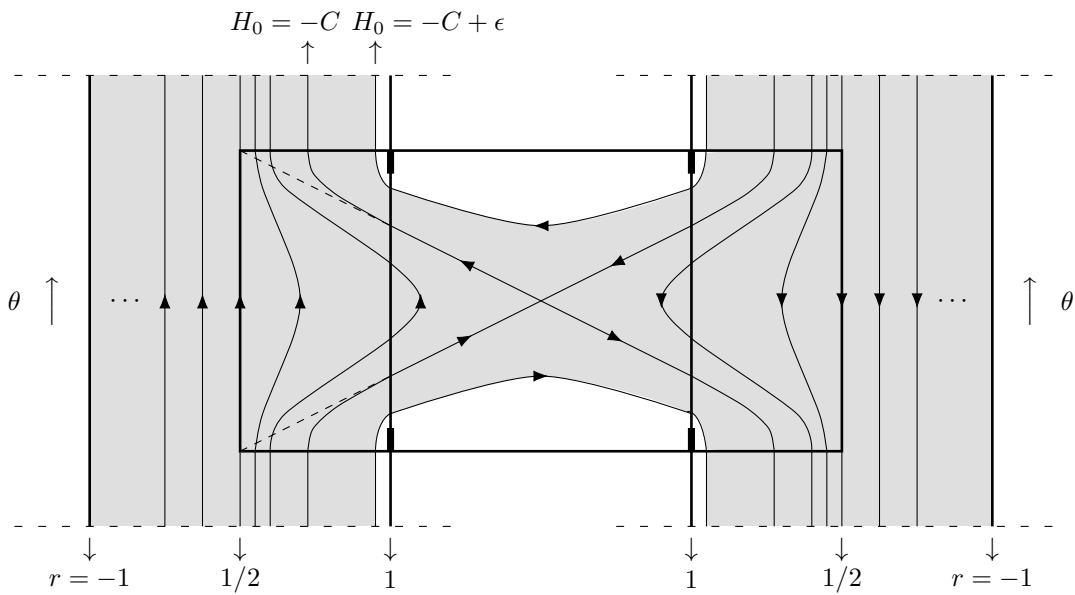
\begin{figure}
	\centering
	\begin{tikzpicture}
	\draw [fill=lightgray!50, draw=none] (-6,-3) rectangle (-2.2,3);
	\draw [fill=lightgray!50, draw=none] (-4,-2) rectangle (4,2);
	\draw [fill=lightgray!50, draw=none] (2.2,-3) rectangle (6,3);
	\draw [fill=white, draw=none] (-1.94,2) ellipse (0.26 and 0.515);
	\draw [fill=white, draw=none] (1.94,2) ellipse (0.26 and 0.515);
	\draw [fill=white, draw=none] (-1.94,-2) ellipse (0.26 and 0.515);
	\draw [fill=white, draw=none] (1.94,-2) ellipse (0.26 and 0.515);
	\draw [fill=white, draw=none] plot [smooth] coordinates {(-2,1.5) (0,1) (2,1.5) (2,2)};
	\draw [fill=white, draw=none] plot [smooth] coordinates {(-2,-1.5) (0,-1) (2,-1.5) (2,-2)};
	\draw [fill=white, draw=none] plot (-2,2) -- (-2,1.48) -- (2.5,2);
	\draw [fill=white, draw=none] plot (-2,-2) -- (-2,-1.48) -- (2.5,-2);
	
	\draw [line width=1pt] (-4,-2) rectangle (4,2);
	\draw [dashed] (-4,-2) -- (-2,-1);
	\draw [dashed] (-4,2) -- (-2,1);
	\draw [line width=1pt] (-2,-3) -- (-2,3);
	\draw [line width=1pt] (-6,-3) -- (-6,3);
	\draw [line width=1pt] (2,-3) -- (2,3);
	\draw [line width=1pt] (6,-3) -- (6,3);
	\draw [loosely dashed] (-7,-3) -- (-1,-3);
	\draw [loosely dashed] (-7,3) -- (-1,3);
	\draw [loosely dashed] (1,-3) -- (7,-3);
	\draw [loosely dashed] (1,3) -- (7,3);
	
	\draw[->-] (-4,-3) -- (-4,3);
	\draw[->-] (-4.5,-3) -- (-4.5,3);
	\draw[->-] (-5,-3) -- (-5,3);
	\draw (-5.5,0) node {$\cdots$};
	\draw[->-] (4,3) -- (4,-3);
	\draw[->-] (4.5,3) -- (4.5,-3);
	\draw[->-] (5,3) -- (5,-3);
	\draw (5.5,0) node {$\cdots$};
	
	\draw (-3.6,-3) -- (-3.6,-2);
	\draw (3.6,-3) -- (3.6,-2);
	\draw (-3.6,2) -- (-3.6,3);
	\draw (3.6,2) -- (3.6,3);
	\draw [line cap=round] plot [smooth] coordinates {(-3.6,2) (-3.4,1.5) (-2,0.5) (-1.6,0) (-2,-0.5) (-3.4,-1.5) (-3.6,-2)};
	\draw [line cap=round] plot [smooth] coordinates {(3.6,-2) (3.4,-1.5) (2,-0.5) (1.6,0) (2,0.5) (3.4,1.5) (3.6,2)};
	
	
	\draw (-3.8,-3) -- (-3.8,-2);
	\draw (-3.8,2) -- (-3.8,3);
	\draw (3.8,-3) -- (3.8,-2);
	\draw (3.8,2) -- (3.8,3);
	\draw [line cap=round] plot [smooth] coordinates {(-3.8,2) (-3.7,1.5) (-3.3,0.5) (-3.2,0) (-3.3,-0.5) (-3.7,-1.5) (-3.8,-2)};
	\draw [line cap=round] plot [smooth] coordinates {(3.8,-2) (3.7,-1.5) (3.3,-0.5) (3.2,0) (3.3,0.5) (3.7,1.5) (3.8,2)};
	
	\draw (-3.1,-3) -- (-3.1,-2);
	\draw (-3.1,2) -- (-3.1,3);
	\draw (3.1,-3) -- (3.1,-2);
	\draw (3.1,2) -- (3.1,3);
	\draw (-2,-1) -- (2,1);
	\draw (-2,1) -- (2,-1);
	\draw [line cap=round] plot [smooth] coordinates {(-3.1,2) (-3.03,1.7) (-2.7,1.4) (-2,1)};
	\draw [line cap=round] plot [smooth] coordinates {(3.1,2) (3.03,1.7) (2.7,1.4) (2,1)};
	\draw [line cap=round] plot [smooth] coordinates {(-3.1,-2) (-3.03,-1.7) (-2.7,-1.4) (-2,-1)};
	\draw [line cap=round] plot [smooth] coordinates {(3.1,-2) (3.03,-1.7) (2.7,-1.4) (2,-1)};
	
	\draw (-2.2,-3) -- (-2.2,-2);
	\draw (-2.2,2) -- (-2.2,3);
	\draw (2.2,-3) -- (2.2,-2);
	\draw (2.2,2) -- (2.2,3);
	
	\draw [line cap=round] (-2.2,2) arc[start angle=180, end angle=255, x radius=0.26, y radius=0.515];
	\draw [line cap=round] plot [smooth] coordinates {(-2,1.5) (0,1) (2,1.5)};
	\draw [line cap=round] (2,1.5) arc[start angle=285, end angle=360, x radius=0.26, y radius=0.515];

	\draw [line cap=round] (-2.2,-2) arc[start angle=180, end angle=105, x radius=0.26, y radius=0.515];
	\draw [line cap=round] (2.2,-2) arc[start angle=15, end angle=90, x radius=0.26, y radius=0.715];
	\draw [line cap=round] plot [smooth] coordinates {(-2,-1.5) (0,-1) (2,-1.5)};
	
	\draw [fill=black] (-2.04,1.7) rectangle (-1.96,2);
	\draw [fill=black] (-2.04,-1.7) rectangle (-1.96,-2);
	\draw [fill=black] (2.04,1.7) rectangle (1.96,2);
	\draw [fill=black] (2.04,-1.7) rectangle (1.96,-2);
	
	\draw (-6,-3.3) node {$\downarrow$}; \draw (-6,-3.7) node {$r=-1$};
	\draw (-4,-3.3) node {$\downarrow$}; \draw (-4,-3.7) node {$1/2$};
	\draw (-2,-3.3) node {$\downarrow$}; \draw (-2,-3.7) node {$1$};
	\draw (6,-3.3) node {$\downarrow$}; \draw (6,-3.7) node {$r=-1$};
	\draw (4,-3.3) node {$\downarrow$}; \draw (4,-3.7) node {$1/2$};
	\draw (2,-3.3) node {$\downarrow$}; \draw (2,-3.7) node {$1$};
	
	\draw (-6.5,0) node {$\Big\uparrow$}; \draw (-7,0) node {$\theta$};
	\draw (6.5,0) node {$\Big\uparrow$}; \draw (7,0) node {$\theta$};

	\draw (-2.2,3.3) node {$\uparrow$}; \draw (-1.5,3.7) node {$H_0=-C+\epsilon$};
	\draw (-3.1,3.3) node {$\uparrow$}; \draw (-3.4,3.7) node {$H_0=-C$};
	
	\draw[->-] [draw=none] (-1.6,-1) -- (-1.6,1);
	\draw[->-] [draw=none] (1.6,1) -- (1.6,-1);
	\draw[->-] [draw=none] (-3.2,-1) -- (-3.2,1);
	\draw[->-] [draw=none] (3.2,1) -- (3.2,-1);
	\draw[->-] [draw=none] (2,1) -- (0,0);
	\draw[->-] [draw=none] (-2,-1) -- (0,0);
	\draw[->-] [draw=none] (0,0) -- (-2,1);
	\draw[->-] [draw=none] (0,0) -- (2,-1);
	
	\draw[->-] [draw=none] (-1,-1) -- (1,-1);
	\draw[->-] [draw=none] (1,1) -- (-1,1);
	
	\end{tikzpicture}
	\caption{The Hamiltonian flow generated by the Morse Hamiltonian (not to scale). The parallel dashed lines glue to form an annulus. The thickened part shows the discontinuities of the cutoff function $\rho$, which we trim back by the level set $H_0^{-1}(-C+\epsilon)$.}
	\label{morse}
\end{figure}
\par
Now we consider the union of both annuli and the handle sets. 
The Hamiltonian is defined as before using cutoff functions. We remark that the right annulus set as depicted in Figure~\ref{morse} has opposite orientation of $r$, so the Hamiltonian flow is also a positive Dehn twist in the annulus set. 
\par
We choose small $\epsilon$ and ``trim'' the domain back to a smooth level set of $H_0^{-1}(-C+\epsilon)$. Since the Hamiltonian $H_0$ in this model only has a critical point in the handle set, the same argument as above shows that we can choose $\epsilon$ such that $H_0^{-1}(-C+\epsilon)$ has smooth boundary, and $H_0$ restricted to this domain is a smooth function. 
As for the return map near the boundary, recall that the size of the Hamiltonian vector field only depends on the slope of the Hamiltonian. Therefore we increase the slope near the boundary such that the return map near the boundary is the identity. 
Since this perturbation will not generate new critical points, the Hamiltonian is still Morse, and generates the return map to be identity near the boundary.
\par
We conclude that the final Hamiltonian $H$ constructed by this process is a Morse Hamiltonian satisfying the properties (1),~(2),~(3) above.

\section{Stability of global surfaces of section}
\label{appendix:stab_ob}
In this appendix, we will prove the following proposition.
\begin{proposition}
\label{prop:stab_gss}
Suppose that $(M,\omega)$ is a symplectic manifold with a Hamiltonian $H$. Assume that the level set $Y=H^{-1}(0)$ is of contact type and admits a global surface of section $\Sigma$ with possibly disconnected binding $B$. Also assume that the periodic orbits of $H$ associated with $B$ are non-degenerate.
Then for a $C^2$-small perturbation $H_\delta$ of $H$, there is an embedded surface $\tilde \Sigma$ that is a global surface of section for the Hamiltonian flow of $H_\delta$.
\end{proposition}

\begin{remark}
Although we have used the words Hamiltonian and contact type to stay in line with the rest of the paper, the statement holds for any $C^1$-small perturbation of a smooth vector field $X$ that admits a global surface of section with non-degenerate binding orbits.
\end{remark}

\begin{proof}
We first outline the argument before going into the details.
Since the global surface of section $\Sigma$ is by definition transverse to the flow on the interior of $\Sigma$, the transversality part is clear for a $C^1$-small perturbation of the vector field generating the flow as long as we stay away from the binding.
We use the linearized flow in order to see that the binding orbits survive a perturbation, and also to see that there is still a strong twist around the binding.

Let us now look at some details.
We consider a perturbation $H+\delta h$ for a $C^2$-small $h$. Except for a neighborhood of the binding, the original open book satisfies the transversality condition for the perturbed Reeb vector field. 
Now consider a small neighborhood $N$ of the binding. Since
the complement of $N$ is a compact region, we can consider $h$ to be small enough that the perturbed Reeb vector field is transverse to the pages of the original open book. Therefore, we only need to consider the inside of the binding neighborhood to show that an open book still exists for the perturbed Hamiltonian. 

We identify the binding neighborhood with $S^1\times D^2$, by introducing coordinates $(x, y, t)$ centered at the original binding. We first look at the Reeb vector field $R_\delta$ for the perturbed Hamiltonian $H_\delta=H+\delta h$. We will actually look at the normalized vector field $\Tilde{R_\delta}$ that has $t$ component identically $1$.
The flow of the normalized vector field is a reparametrization of the flow of the original vector field. 
We now apply an implicit function theorem argument, that shows for small $\delta$, there exists a periodic orbit $\gamma_\delta$ near the original periodic orbit. 

Define $X(x,\delta)$ to be the Hamiltonian vector field of the perturbed Hamiltonian $H+\delta h$. 
We can take a point $p$ on the binding. 
Take a local surface of section $S$ transverse to the binding at $p$. We define the return map of the Hamiltonian vector field at perturbation $\delta$ to be $\varphi(x,\delta)$. 
Since the binding orbit is non-degenerate, the linearization of $\varphi(x,\delta)$ at $(p,0)$ does not have $1$ as an eigenvalue. 
Now consider the map $S\times[0,\delta_1]\to S$, where $[0,\delta_1]$ is the domain of the parameter $\delta$. 
From the discussion above, we have that the matrix $\frac{\partial}{\partial \delta}\varphi-I$ is non-singular at $(p,0)$. 
Therefore from the implicit function theorem, we can find a point $\gamma_\delta(0)$ in $S$ such that $\varphi(\gamma_\delta(0),\delta)=\gamma_\delta(0)$ for small $\delta$. 
This implies that the orbit $\gamma_\delta(t)$ through this point is a smooth periodic orbit for the perturbed Hamiltonian $H+\delta h$.

We denote the coordinates of this orbit by $(\gamma_{\delta,x}, \gamma_{\delta,y},\gamma_{\delta,t})$. 
Then we introduce new variables $\tilde u, \tilde v$ by putting $\tilde u=x-\gamma_{\delta,x}(\gamma_{\delta,t}^{-1}), \tilde v=y-\gamma_{\delta,y}(\gamma_{\delta,t}^{-1})$. 
This $u,v$ measures the distance from the perturbed orbit $\gamma_\delta$, and $z$ is a reparametrization of the $t$ coordinate. 
These coordinates $(\tilde u,\tilde v,z)$ are defined in a tubular neighborhood $\nu_Y(\gamma_0)$ of the unperturbed periodic orbit $\gamma_0$, but their behavior on the boundary of $\nu_Y(\gamma_0)$ depends on $\delta$.
To fix this, we use a cutoff function $\rho$, which equals $1$ on $\gamma_\delta$ and vanishes on a neighborhood of $\partial \nu_Y(\gamma_0)$.
Then we put $u_\delta=x-\rho\gamma_{\delta,x}(z), v_\delta=y-\rho\gamma_{\delta,y}(z)$.

Define the modified ``open book'' projection $\theta_\delta=\frac{(u_\delta, v_\delta)}{\sqrt{u_\delta^2+v_\delta^2}}\in S^1$. 
This map coincides with the original open book projection near $\partial \nu_Y(\gamma_0)$ and can hence be extended smoothly to $Y$ by using the original open book projection.
By cutting out the global surface of section, we can lift this map to $\R$, which we will do to have a convenient description of the derivative; we will continue to write $\theta_\delta$, also for this lifted map.
Since the unperturbed Reeb vector field $R_0$ is transverse to the interior of the global surface of section, we can find $C>0$ such that $R_0(\theta_0)>2C$ on $Y \setminus \nu_Y(\gamma_0)$ (away from all binding orbits).
As $R_\delta$ is $C^1$-close to $R_0$, we still have $R_\delta(\theta_\delta)>C$ for $0<\delta<\delta_1$ if we choose $\delta_1$ sufficiently small.

It hence suffices that to show that we have transversality on a neighborhood of the binding orbit, $\nu_Y(\gamma_0)$.
To analyze this, consider the smooth $1$-form 
$$
\Omega_\delta=u_\delta dv_\delta-v_\delta du_\delta.
$$
We observe that 
$$
d \theta_\delta =\frac{\Omega_\delta}{u_\delta^2 +v_\delta^2},
$$
so $R_\delta(\theta_\delta)>0$ is equivalent to $\Omega_\delta(R_\delta)>0$.

Since $\Omega_\delta (R_\delta)$ is a smooth function of $p=(u,v,t)$ and $\delta$, we consider a Taylor expansion in the $u,v$-coordinates and $\delta$. 
The $0$-th order term in this expansion is
\begin{equation*}
\Omega_0(R_0)(u,v,t) = C_t(u^2+v^2)+o(u^2 + v^2).
\end{equation*}
This can be seen most easily from the explicit form of the Reeb vector field~\eqref{eq:Reeb_vf}, but below we shall see that such an expression follows for all small $\delta$ by analyzing the linearized flow. 
We make the following two claims.~\\
\noindent
{\bf Claim 1: } $\Omega_0(R_0)\geq 0$ and $\Omega_0(R_0)$ vanishes only along $\gamma_0$. 
 ~\\
\noindent
{\bf Claim 2: }
there is a uniform (i.e.~independent of $\delta$) neighborhood $N$ of $\gamma_0$, and $\delta_2$ such that for $\delta \in [0,\delta_2]$ the following hold.
\begin{itemize}
\item $\gamma_\delta\subset N$
\item for all $p\in N$ we have $\Omega_\delta(R_\delta)(p) \geq 0$ and $\Omega_\delta(R_\delta)(p)=0$ if and only if $p$ in the image of $\gamma_\delta$.
\end{itemize}
The first claim is clear. 
We verify the second claim by analyzing the linearized flow.

Put $P=(U,V,Z)$, and set $p =\gamma_\delta+\epsilon P$.
The flow equation for $p$ is 
$$
\frac{d p}{dt} =\tilde R_\delta(p),
$$
and by expanding in $\epsilon$ we obtain the linearized equation $\frac{d P}{dt} =\epsilon \nabla_P \tilde R_\delta +o(\epsilon)$.
Since we only need the component normal to $\gamma_\delta$, we will use the following matrix representation for the normal component of the linearized flow.
$$
\left(
\begin{array}{c}
\dot U \\
\dot V
\end{array}
\right)
=
A_\delta \left(
\begin{array}{c}
 U \\
 V
\end{array}
\right)
,
$$
where $A_\delta$ is a time-dependent matrix.  
We now compute the value of $\Omega_\delta(\tilde R_\delta)$ at $p$ using the above definitions.
\begin{align*}
\Omega_\delta(R_\delta) (p)&=
(u_\delta dv_\delta-v_\delta du_\delta)(R_\delta)\\ 
&= \epsilon^2(U \dot V - V \dot U) + o(\epsilon^2)\\
&=
\left(\begin{smallmatrix}
\epsilon U \\
\epsilon V
\end{smallmatrix}\right)^t
\left(\begin{smallmatrix}
0 & 1\\
-1 & 0
\end{smallmatrix}\right) A_{\delta}
\left(\begin{smallmatrix}
\epsilon U \\
\epsilon V
\end{smallmatrix}\right)+o(\epsilon^2).
\end{align*}
For fixed $t$, we know from the above that $\left(\begin{smallmatrix}
0 & 1\\
-1 & 0
\end{smallmatrix}\right)A_0$ is positive definite, so $\left(\begin{smallmatrix}
0 & 1\\
-1 & 0
\end{smallmatrix}\right) A_\delta$ is too, for sufficiently small $\delta$.
This settles the second claim.

To complete the proof, we argue by contradiction.
Suppose that for all $\delta>0$, there is a point $p_\delta \notin \gamma_\delta$ such that $\Omega_\delta(\tilde R_\delta)(p_\delta)=0$.
We obtain a sequence $\delta_n$ converging to $0$ and, by compactness, a converging sequence $p_n$, such that $\Omega_{\delta_n}(\tilde R_{\delta_n})(p_n)=0$ for all $n$.
By Claim 1, we see that $p_\infty$ is in the image of $\gamma_0$.
But this means that $p_n$ lies in $N$ for sufficiently large $n$, contradicting Claim 2.
This completes the proof.
\end{proof}

\end{appendices}

\end{document}